
\input amstex

\input psfig
\documentstyle{amsppt}
\baselineskip=24pt
\hsize=6.5 true in 
\vsize=8.75 true in
\rightheadtext{Ratio geometry, rigidity and the scenery process}
\widestnumber\key{deM-vS}
 
\define\r{{\Bbb R}}

\define\wt{\widetilde}
\define\wh{\widehat}
 
\topmatter
\title
Ratio geometry, rigidity and the scenery process for hyperbolic
Cantor sets.
\endtitle
\author
Tim Bedford
and
 Albert M. Fisher
\endauthor
\affil
Delft University of Technology;
SUNY at Stony Brook
\endaffil
\date{June 26, 1994}
\enddate
\abstract{
 Given a $\Cal C^{1+\gamma}$ hyperbolic Cantor set
$C$, we study the sequence $C_{n,x}$ of Cantor subsets which nest down
toward a point $x$ in $C$.  We show that $C_{n,x}$ is asymptotically
equal to an ergodic Cantor set valued process.  The values of this
process, called \it limit sets\rm , are indexed by a H\"older
continuous set-valued function defined on D. Sullivan's dual Cantor
set. We
show the limit sets are themselves 
$\Cal C^{k+\gamma},~\Cal C^\infty$ or $\Cal C^\omega$
hyperbolic Cantor sets, with the highest degree of smoothness which
occurs in the $\Cal C^{1+\gamma}$ conjugacy class of $C$.  The proof of this leads
to the following rigidity theorem:  if two $\Cal C^{k+\gamma},~\Cal C^\infty$
or $\Cal C^\omega$ hyperbolic Cantor sets are $\Cal C^1$-conjugate, then the
conjugacy (with a different extension) is in fact already
$\Cal C^{k+\gamma},~\Cal C^\infty$ or $\Cal C^\omega$.
}
Within one $\Cal C^{1+\gamma}$
conjugacy class, each smoothness class is a Banach manifold, which is
acted on by the semigroup given by rescaling subintervals. Conjugacy
classes nest down, and contained in the intersection of them all 
is a compact set which is the attractor for the semigroup:
 the collection of limit sets. 
Convergence is exponentially fast, in the $\Cal C^1$ norm.
\endabstract
\endtopmatter
 
\document
\NoBlackBoxes

\heading {\bf Introduction.} \endheading
Consider
 the sequence $C_{n,x}$ of Cantor subsets which nest down
toward a point $x$ in a hyperbolic Cantor set $C\subseteq [0,1]$, and which
have
been affinely rescaled to have left and right endpoints at $0$ and
$1$.  We wish to describe how the geometry of these sets changes as $n$
increases. If $C$ is a linear set like the middle-third set this is not so
difficult to do (we always get just another copy of $C$!) but, as we shall
see, with nonlinearity the behavior of this ``scenery process'' gets
much more interesting.
 
A different way to describe the small-scale structure of $C$ is by
the \it scaling function\rm , introduced by Feigenbaum for a specific
class of examples, and studied by D. Sullivan in the present setting
of $\Cal C^{1+\gamma}$ hyperbolic Cantor sets.

A third approach is to define a flow, the continuous dynamics of
which reflect the geometrical notion of zooming continuously down
toward a point.  Ergodicity of the flow implies analogues of the
Lebesgue density theorem, proved in \cite{BF 1} for Brownian zero sets
and hyperbolic $\Cal C^{1+\gamma}$ Cantor sets; the number one gets
(order-two density) is a conformal
invariant and provides a measure of
the \it lacunarity\rm~ of the fractal; compare  \cite{Mand}.
This \it scenery flow\rm
~is constructed for hyperbolic Cantor sets in \cite{BF 2}, \cite{BF 3},
 for hyperbolic
 Julia sets in \cite{BFU}, and
for limit sets of 
geometrically finite Fuchsian and Kleinian groups in \cite{Fi 2}.
The associated translation (rather than dilation) 
scenery flows are studied for 
hyperbolic $\Cal C^{1+\gamma}$ Cantor sets
and the Fuchsian limit sets in 
 \cite{Fi 2} and \cite{Bu-F}, where theorems like those of \cite{Fi 1} 
are proved 
(an order-two ergodic theorem, and infinite-measure unique ergodicity).

In the present paper we will take a viewpoint close to that of Sullivan in
\cite{Su 1}, constructing the scenery process in a similar way to the
scaling function.  It is also possible to work in the other
direction; in the later papers \cite{BF 2}, \cite{BF 3}
we show how to derive the
scaling function from the scenery flow, and conversely how to
construct the scenery flow from the scaling function.  From this
point of view the scenery process will be seen as an intermediate
object, serving to connect the scenery flow with the scaling function.

Sullivan's main motivation in \cite{Su 1}
was to begin to develop a Teichm\"uller theory for Cantor sets, for use in
a new, more ``conceptual'' proof of 
the Feigenbaum-Coullet-Tresser conjectures (see \cite{Su 2} and \cite{deM-vS}). 
(Lanford's proof \cite{L}
 uses (rigorous) computer-assisted estimates). 
A first step would be to classify 
differentiable structures on the attracting Cantor sets of folding maps. 
This classification should be in terms of some invariant which would
serve as a ``modulus'' of the structure; the next step would be to 
put a complex structure on the set of moduli as in classical 
Teichm\"uller theory, see \cite{Su 3}. Now for the particular case of the folding map which is 
the
Feigenbaum-Coullet-Tresser
renormalization fixed point, the Cantor set not only has a folding 
dynamics but also a hyperbolic $\Cal C^{1+\alpha}$
dynamics. (This observation is attributed by Rand \cite{Ra}
to Misiurewicz). One is thus led to the following separate question,
which is the subject of \S \S 1-3 of \cite{Su 1}: for general hyperbolic
$\Cal C^{1+\alpha}$ Cantor sets, can one classify differentiable
structures? Sullivan shows this can be done, with the
``modulus" being a bounded H\"older scaling function.

Our own main focus is somewhat different. We want to describe
the \it exact \rm geometry, at small  scales, of the Cantor sets,
 whereas to the  differentiable structure, all
smoothly equivalent 
Cantor sets will look the same. 
However for this purpose also, the scaling function
contains precisely the information one needs.

Our main theorems (Theorems 5.4,  7.4 and 7.5)
concern respectively the scenery process, the smoothness of limit sets, 
and $\Cal C^{k+\gamma}$ rigidity.
We summarize the totality of the resulting picture.
Given one $\Cal C^{1+\gamma}$ hyperbolic Cantor set, consider
the collection of all Cantor sets which are $\Cal C^{1+\gamma}$ 
conjugate to it.
Within this collection is a distinguished subcollection, its limit sets.
The free semigroup on two generators acts on the conjugacy class (by 
rescaling subsets of the next level); the 
limit sets
are an attractor for this action, 
and the scenery process can be described as what one sees when 
walking out a branch
of the tree of the semigroup.
Limit sets are exactly the \it ratio sets \rm (see \S 2) 
built from the associated scaling
function. Within the big collection are subcollections with higher
degrees of smoothness. 
The big collection forms an infinite-dimensional Banach manifold,
naturally identified with a factor of
the $\Cal C^{1+\gamma}$-diffeomorphisms of the 
interval,
after one Cantor set has been chosen as a base point. (The diffeomorphisms 
are the conjugacies to this set).  
The subcollections  
nest 
down as smoothness increases, and
 by rigidity these smoothness classes are conjugacy classes 
as well. Contained in the intersection of them all is
the collection of  limit sets, with
the highest possible smoothness.
Choosing one of them as a common base point,  
 these  subcollections  are naturally identified
with the $\Cal C^{k+\gamma}$ diffeomorphisms of the interval. Each is a
Banach manifold in its own topology, and is dense in a larger collection 
with respect to its topology. The free semigroup acts on each manifold. Its 
 points are drawn exponentially fast in the $\Cal C^1$ norm toward 
the common attractor: the collection of limit sets, which
form a compact subset of the Banach manifold.

In the course of our paper
we give careful proofs of several of Sullivan's theorems (\cite{Su 1} 
is extremely sketchy).
In some cases our different point of view leads us to different arguments 
from those indicated in 
\cite{Su 1}. We will describe our approach and results 
more fully after a further 
explanation of Sullivan's ideas.

\subhead
{Sullivan's  differentiable structures}
\endsubhead

We begin with an ordered topological Cantor set,
i.e\. a space which is homeomorphic and order
 isomorphic to the usual middle--third Cantor set. For convenience we use 
$\sum^+\equiv\Pi_0^{\infty}\{0,1\}$,
together with the product topology, and with the 
 lexicographic order.
{\bf Charts}
are defined to  be order--preserving homeomorphisms into $\r$; 
two charts $\zeta,\ 
\xi$ are 
$\Cal C^{k+\gamma}$
{\bf compatible} if $\zeta\circ\xi^{-1}$
extends, with that degree of smoothness, to a diffeomorphism defined on 
neighborhoods of the embedded sets. A 
{\bf linear} $\Cal C(k,\alpha)$
{\bf differentiable structure} on $\Sigma^+$
will be a {\bf maximal atlas} (a maximal compatible collection of charts).
Here, following \cite{Su 1},  
$\Cal C(k,\alpha)$ 
denotes all maps which are
$\Cal C^{k+\gamma}$
for some $\gamma\in (0,1]$.
Therefore, a 
$\Cal C(k,\alpha)$ 
linear differentiable structure determines and is determined by 
a class of Cantor sets
embedded in the real line, equivalent by $\Cal C(k,\alpha)$ 
 changes of coordinates.

For simplicity, we are restricting our attention to charts which are order 
preserving and globally   
defined. We mention that the word ``linear'' is being used 
in two ways: when dynamics is introduced on these sets, it will
 usually be nonlinear; the differentiable stuctures are called linear
because they come from embeddings in the line.
(Alternative theories might have charts mapping 
$\Sigma^+$ to a product of Cantor
sets, or to a subset of some fractal curve!)

Via the homeomorphism from 
$\Sigma^+$,
an embedded set $C$ comes equipped with  the  dynamics of
 the shift map 
$\sigma$ on $\Sigma^+$. The set also inherits from $\Sigma^+$
a
nested hierarchy of intervals, corresponding to finite words (\it cylinder
sets\rm) in $\Sigma^+$.
Sullivan uses the shift map to define 
$\Cal C(1,\alpha)$
\it hyperbolic Cantor sets \rm
(see \S 1 below) 
and the nested intervals to define the \it ratio geometry \rm of a Cantor set.
This assigns
to each interval the triple $(l,g,r)$ of length ratios of the left
subinterval, middle gap and right subinterval respectively.
The hypothesis that an embedded Cantor set $C$ is 
hyperbolic $\Cal C(1,\alpha)$
is enough to show that a limiting ratio geometry, recorded by the 
\it scaling function, \rm exists. 
The ratio geometry is bounded away from $0$ and $1$, a condition called 
\it bounded geometry\rm; hence the limiting values $(a,b,c)$ 
are also in the interior of 
the simplex $\Delta=\{(a,b,c):\ a+b+c=1\}$. 
Convergence to the scaling 
function is taken along  inverse branches of $\sigma$, which are
indexed by points of an abstract topological Cantor set called 
the \it dual Cantor set\rm. Thus the scaling function maps the dual
Cantor set to $\Delta$, onto a compact subset of its interior.
Convergence is exponentially fast, and the scaling function 
is  H\"older continuous. 

Locations in the Cantor set correspond to forward images under $\sigma$,
since the digits of $\Sigma^+$ tell whether the orbit of a point lies 
in the left or right third of $C$. To study smoothness of a conjugacy
or an expanding map, one expects of course to use the locations to
estimate difference quotients. 
However since convergence to the scaling 
function is taken along {\it inverse} branches of $\sigma$,
as the scale 
gets smaller and smaller,  the locations jump all over the set.

The first remarkable result from \S\S $1-3$ of \cite{Su 1} is that while
indeed one cannot compute the derivative of the shift map from the scaling
function, nevertheless this function contains complete information about
$\Cal C(1,\alpha)$ 
differentiable structures.

More precisely, one has the following. As we have already mentioned,
(1) a hyperbolic $\Cal C(1,\alpha)$ Cantor set has a bounded H\"older
scaling function. Next (2) this depends only on the differentiable
structure i.e\. it is the same for $\Cal C(1,\alpha)$
conjugate Cantor sets. Conversely (3) an embedded Cantor set which
has a bounded, H\"older scaling function is in fact $\Cal C(1,\alpha)$
hyperbolic. Finally, (4) in this case the $\Cal C(1,\alpha)$ differentiable
structure is determined by the scaling function. In other words two
hyperbolic $\Cal C(1,\alpha)$ Cantor sets with the same scaling function are
$\Cal C(1,\alpha)$ conjugate. In summary, the bounded H\"older scaling
function gives an intrinsic characterization  of the  differentiable 
structure, in the sense that no embedding need be specified.

Now furthermore, quoting \cite{Su 1}: ``\dots if the structure admits
a $\Cal C(k,\alpha)$ refinement so that the shift is $\Cal C(k,\alpha)$,
this structure is also determined uniquely by the same scaling function
\dots". Stated as a result about representatives instead of the entire
equivalence class, this can be interpreted as a rigidity theorem: if
two $\Cal C(k,\alpha)$ hyperbolic Cantor sets are conjugate by a map
which is $\Cal C(1,\alpha)$, then that map (possibly with a different 
extension to the gaps) is in fact already $\Cal C(k,\alpha)$.

\subhead
{Summary of results}
\endsubhead

We include in this paper careful statements and proofs 
in particular of (1), (2), (4) above, and of
$\Cal C(k,\alpha)$  
rigidity. 
(We mention that our use of the term ``rigidity"
is different  from that in \S $5$ of \cite{Su 1}).
Since we are interested in the geometry of
representatives rather than the equivalence class, all these results 
are stated in terms of the 
conjugacy of embedded sets rather than the classification of
differentiable 
structures. As we said above, the reason for this emphasis 
is that our primary goal is to study 
the scenery process, and all the sets in the scenery process are the same up
to conjugacy.

A small technical difference to \cite{Su 1}
is that we use $\Cal C^{k+\gamma}$ rather 
than $\Cal C(k,\gamma)$ throughout.
We do this because it  
gives sharper statements. 
Thus e.g\. for rigidity we show
that
$\Cal C^1$
conjugacy 
implies $\Cal C^{k+\gamma}$
conjugacy.

In part because of our change in focus, we give a different proof from that 
suggested in \cite{Su 1} of (4). Each approach has its own advantages. Sullivan's
method, a direct estimate of the derivative by difference quotients using
sums of gap lengths, gives a unified way of proving (3) as well as (4).
However one also needs to cite an extension lemma, which is not included in 
\cite{Su 1}. On the other hand our approach gives a unified treatment of 
$\Cal C^1$
conjugacy and rigidity, and avoids calling on the separate extension lemma.
We do not prove (3) here, but will give a full proof (along the lines
of \cite{Su 1}) elsewhere.

Our proof of rigidity is intimately connected to the study of the 
scenery process. 
We proceed as follows.

First
we use limiting conjugacies to construct a set-valued analogue of the
scaling function.  This function, $y\mapsto C^y$, is defined for $y$ in the
dual Cantor set, is H\"older continuous with respect to a metric derived from
the corresponding Hausdorff measures and has as its range a
compact subset of the collection of all subsets of $[0,1]$ in that
measure metric, and also in the
Hausdorff metric on sets.  The scaling dynamics enters by interpreting
the dual Cantor set as the past of the natural extension of the
expanding map on $C$; the scenery process $C_{n,x}$ is then
asymptotically given by evaluating the shift on any extension $\underline
x =(y,x)$ of $x\in C$.  Since the limit sets were constructed by
conjugacies, one can apply a lemma from the appendix of \cite{Su 1}
to help determine their degree of smoothness:  we show they have the highest
degree of smoothness ($\Cal C^{k+\gamma}$ for some $k\geq 1$, $\Cal C^\infty$
or $\Cal C^\omega$) which occurs in the $\Cal C^{1+\gamma}$ conjugacy class of $C$.

Next, the proof of rigidity follows as a corollary. Given two hyperbolic
$\Cal C^{k+\gamma}$
Cantor sets, if they are $\Cal C^1$
conjugate they have the same scaling function. Hence they have the same 
limit sets,
which are \it ratio Cantor sets \rm constructed from this function. Choosing
one of these to act as an intermediary, the composition should also be
$\Cal C^{k+\gamma}$. 
However the maps may be defined differently on the gaps, which would lead back 
to the extension problem mentioned before.
But now one has a simpler solution: a choice is made on the middle third,
and the rest of the definition follows automatically from the dynamics. 
This completes the proof of rigidity. 
In summary: if
two $\Cal C^{k+\gamma},~\Cal C^\infty$ or $\Cal C^\omega$ hyperbolic Cantor 
sets are $\Cal C^1$
conjugate, then this conjugacy (with a different extension) is
already $\Cal C^{k+\gamma},~\Cal C^\infty$ or $\Cal C^\omega$ respectively.

This leads, then,  to the overall picture which is
summarized at the end of the first part
of the Introduction.

Now  one knows, from the rigidity theorem, 
that  the maximum degree of smoothness  occuring in the 
$\Cal C^{1+\gamma}$ conjugacy class
should be encoded somehow in the scaling function.
Work of Tangerman and Przytycki
gives one way of recovering that information \cite{T-P}.  
 A. Pinto and
 D. Rand (\cite{P-R 2}, \S 5 and personal communication) and   Dennis 
Sullivan
(personal communication) have suggested other approaches, in a 
related situation.
It would be nice to understand in a unified way
these different
points of view.

The rigidity theorem is also  stated by
Tangerman and Przytycki; it is proved  
as a corollary of their 
main result. 
Their approach  is quite different from ours and in particular does 
make use of 
Whitney's 
Extension Theorem.  (We became aware of their preprint after  
the first version of this paper - an IHES preprint, July 1992 -
 was  completed).

 Rand introduces the notion of  a 
 \it  
Markov family \rm to help
study the relationship between scaling functions and smooth conjugacy
in  situations where one has a \it sequence \rm of expanding 
maps, 
rather than a single map.
The examples studied in \cite{Ra}
, \cite{P}, and \cite{P-R 1,\ 2} include 
certain circle diffeomorphisms and folding maps. 
See \cite{AF} for some related developments.

Interesting work on the  small-scale geometry 
of certain fractal sets, in quite different settings,
has been done by
Hillel Furstenberg and Tan Lei. 
Tan Lei in \cite{T} proves the beautiful theorem that certain nonhyperbolic
Julia sets, corresponding to Misiurewicz points in the boundary of the 
Mandelbrot set $\partial M$, are  
asymptotically the same as $\partial M$ at that point.  These points 
form a countable dense subset  of 
$\partial M$, yet the general case is still far from completely understood.
This type of asymptotic limit, as well as what we have called here
limit sets, provide examples of Furstenberg's   general notion of the
{\it microsets} of a subset of Euclidean space (lectures 
and personal communication). 
These are by definition   all the
limiting sets given by rescaling nested subsets by a 
sequence of  affine expansions.
Furstenberg applies this in a continuation of the analysis
 begun in \cite{Fu} for determining the 
Hausdorff dimension of
certain sets: intersections of generic translates of linear Cantor sets, and  
intersections of linear 
Cantor sets in the plane with foliations of straight lines at a generic slope.
His study of these matters is related to
the ``times $2$ times $3$'' circle of 
problems in Ergodic Theory. An interesting and important area of research
is to develop similar results in a nonlinear setting, e.g. for general
smooth foliations or for nonlinear Cantor sets.

\subhead{Acknowledgements}
\endsubhead
We wish to thank our colleagues and friends
 for conversations, encouragement and inspiration regarding this 
and  related projects. 
We
 give special thanks to
M. Urbanski, B. Mandelbrot, S. Kakutani, 
and D. Sullivan. 
 The second-named
author would also like to thank the
Dynamical Systems seminar at Memphis
State University for their encouragement to give a series of lectures on
these topics (Spring 1990), and Yale University, MSRI, IHES,  CUNY,  the CNRS,
 Universit\'e Paris-Nord, 
and SUNY at Stony Brook for their support while the
paper was being written.
 

\subhead{\S 1 Two ways of building Cantor sets}
\endsubhead
\subhead{\S 1.1 Hyperbolic Cantor sets, Hausdorff and Gibbs measures}
\endsubhead
We start with the usual middle-third Cantor set.
Let $S$ denote
the 2-1 map on the middle-third set $C$ defined by $x\mapsto
3x(\text{mod} 1)$.  The Hausdorff dimension of $C$ is $d=\log 2/\log 3$;
writing $H^d$ for $d$-dimensional Hausdorff measure and $\mu$ for
the restriction $\mu=H^d|_C$, we recall that $\mu$ is a Borel probability
measure (total mass =1) which is invariant under $S$.  The triple
$(C,S,\mu)$ is canonically isomorphic to the one-sided Bernoulli
left 
shift $\sigma$ on $\Sigma^+\equiv\underset 0\to{\overset \infty \to \prod}
\{0, 1\}$, with infinite $(\frac 12, \frac 12)$ coin-tossing
measure; the correspondence is given by $\pi:(x_0 x_1\dots)\mapsto x $
where $x\in C$ has ternary expansion
$$x=\sum_{i=1}^\infty 2x_i 3^{-i}.
$$
A {\bf hyperbolic $\Cal C^{1+\gamma}$ Cantor set} $C$
by definition also has an expanding dynamics
$S:C\longrightarrow C$, but now instead of having straight lines as for
$3x(\text{mod} 1)$, the graph of $S$ may be nonlinear:

 
\


\centerline{\psfig{figure=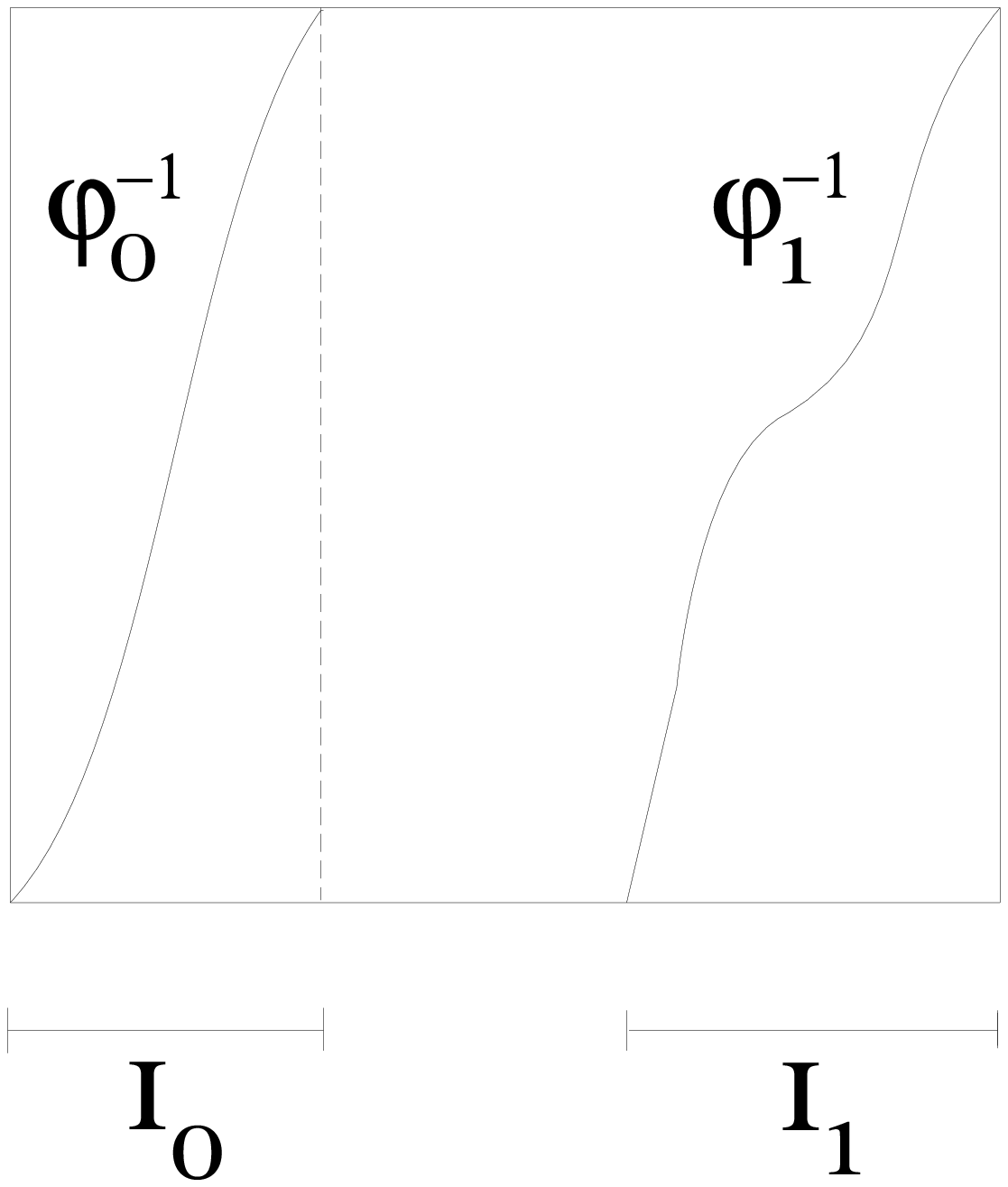,height=1.9in
}}

\smallskip\smallskip
\centerline{[Figure 1]}

\

To construct such a set, one defines $C$ as a limit from two contraction 
mappings
$\varphi_0, \varphi_1:I\longrightarrow I$.
We consider first the case where these maps
are orientation-preserving, and are \it strict\rm~ contractions in the sense 
that
the derivatives satisfy
$0<\alpha<D\varphi_i<\beta<1$. We also require that
$$0=\varphi_0(0)<\varphi_0(1)<\varphi_1(0)<\varphi_1(1)=1.$$

This implies that the
intervals $I_0\equiv\varphi_0(I), I_1\equiv\varphi_1(I)$ are disjoint.
We assume that $\varphi_0, \varphi_1$ are
$\Cal C^{1+\gamma}$ maps for some $\gamma\in(0,1]$.
Here $\Cal C^{k+\gamma}$
means the $k^{\text{th}}$ derivative $D^k \varphi_i$ is H\"older continuous 
with
exponent $\gamma$; note that $\Cal C^{1+1}$ means $D\varphi_i$ is Lipschitz,
so $\Cal C^2$ implies $\Cal C^{1+1}$ (by compactness) but not conversely.
(Exponent $\gamma>1$ is excluded because in that case, since the domain $I$
is connected, $\varphi_i$ is identically constant hence immediately of order
$\Cal C^\infty$ - while the whole purpose of H\"older conditions is to have
\it intermediate\rm~ grades of smoothness).
Our convention for H\"older continuity will be: if we are given that 
$f$ satisfies 
$|f(x)-f(y)|\leq
c_0|x-y|^\gamma$,
then we say $f$ is {\bf H\"older continuous} with {\bf H\"older constant} 
$c_0$
and 
{\bf H\"older exponent} $\gamma$.
We remark that Sullivan in \cite{Su 1} works instead with
$\Cal C(k,\alpha)$ maps; these are defined to be the collection of
$\Cal C^{k+\gamma}$ maps for all $\gamma\in (0,1]$.

We define $S:I_0\cup I_1\to I$ to be the map with inverse
branches $\varphi_0, \varphi_1$. Note that since $D\varphi_i$ are bounded
away from $0$ and $\infty$, it follows that $S$ is $\Cal C^{1+\gamma}$
with same H\"older exponent, but with different H\"older constant.
 
Inductively, form
$$
I_{x_{0}\dots x_n} =
\varphi_{x_{0}}(\varphi_{x_{1}}\dots(\varphi_{x_{n}}(I)))$$
where $x_k\in \{0, 1\}; \; \bigcup I_{x_0\dots x_n}$ (union over all
choices, with $n$ fixed)
is the $n^{\text{th}}$
level approximation  to the Cantor set, $C$, defined as
$$
C=\bigcap_{n=0}^\infty \bigcup I_{x_0\dots x_n}.
$$
The restriction of the map $S$ to $C$ maps $C$
to itself and is (just as for the middle-third set) conjugate to
the Bernoulli shift $(\Sigma^+,\sigma)$, via the map 
$\pi:(x_0x_1\dots)\mapsto x$,
where $x$ is the unique element of $\cap_{n=0}^{\infty}I_{x_0\dots x_n}$.
 
A set $C$ together with map $S:I_0\cup I_1\to I$,
defined in this way from strict contractions
$\varphi_0,\varphi_1$,
will be called a {\bf strictly hyperbolic $\Cal C^{1+\gamma}$
Cantor set (with map)}. 

Sometimes we are only interested in the dynamics on $C$ itself.
 Knowledge of this
{\bf restricted map} $S|_C$ is equivalent to knowing how  $C$ is coded 
by $\Sigma^+$.
We will refer to a set together 
with this labeling as a {\bf marked} Cantor set.
When we forget about this coding, $C$ will be referred to as
the {\bf underlying} Cantor set of $(C,S)$.

More generally, a {\bf hyperbolic } $\Cal C^{1+\gamma}$ Cantor set is
defined as follows. Again we assume that
$\varphi_0,\varphi_1: I\to I$ are order-preserving $\Cal C^{1+\gamma}$
diffeomorphisms such that
$\varphi_0(0)=0,
\varphi_1(1)=1$, and $\varphi_0(1)<\varphi_1(0)$.
We also assume as before that there exists $ \alpha $ with 
$0<\alpha<D\varphi_i<1.$ 
However now the upper bound is replaced by one of the two equivalent conditions
which follow. We write $\varphi_{x_0\dots x_n}
\equiv \varphi_{x_0}\circ\cdots \circ
\varphi_{x_n}$.

\proclaim{Proposition 1.1.1}
The following are equivalent:
 
\item{(i)}
$\exists \beta<1$ and $N\geq 0$ such that for all $n>N$, for any 
$x_0\dots x_n$,
$
D\varphi_{x_0\dots x_n} <\beta^n
$
\item{(ii)}
$\exists c>0 $ and $\tilde\beta<1$ such that for all $n$, 
for any $x_0\dots x_n$,
$
D\varphi_{x_0\dots x_n} <c\tilde\beta^n.
$
\endproclaim
\demo{Proof}
To pass from (ii) to (i), any $\beta$ will work which satisfies
$\tilde\beta<\beta<1$.
For the converse, take 
$c=\text{ max}\{(D\varphi_{x_0\dots x_N})\cdot\beta^{-N}\}$
and $\tilde\beta=\beta$.
\ \ \qed\enddemo
 
Two $\Cal C^{1+\gamma} $ hyperbolic Cantor sets
$C,\ \widetilde C$  with maps $S, \ \widetilde S$
will be said to be {\bf continuously conjugate}, or 
$\Cal C^0-${\bf conjugate},  if there is an 
order-preserving homeomorphism $\Phi :I\to I$ such that for all
$x\in I_0\cup I_1$,
$\widetilde S\circ \Phi(x)=\Phi\circ S(x)$.
We say
$C,\ \widetilde C$
are {\bf $\Cal C^{k+\gamma}, \ \Cal C^\infty, \ \Cal C^\omega$ conjugate} if
$\Phi$ and its inverse have that degree of smoothness.
Note that from the definition, for each $n$,
$\Phi(I_{x_0\dots x_n})=\widetilde I_{x_0\dots x_n}$.
Therefore
the conjugacy induces the identity map on the corresponding shift spaces.
 
When $\Phi$ is  defined (as above) on all of $I$  we will also call it a
{\bf full} conjugacy.   
A {\bf restricted}  conjugacy  is a conjugacy between the Cantor sets 
which can be extended to a full conjugacy.

\proclaim{Lemma 1.1.2} Let $(C,S)$ be a hyperbolic Cantor set, and
assume $\Phi:I\to I$ is a $\Cal C^{1+\gamma}$ diffeomorphism. Define sets
$\wt C\equiv\Phi(C)$, $\wt I_i\equiv\Phi(I_i)$ for $i=0,1$,
and define the map $\wt S:\wt I_0\cup\wt I_1\to I$ by
$
\wt S\equiv
\Phi\circ S\circ\Phi^{-1}
.$
Then $(\wt C,\wt S)$ is also a hyperbolic $\Cal C^{1+\gamma}$ Cantor set.
\endproclaim
\demo{Proof}
This is immediate. Note that using condition (i), $\tilde\beta$ stays
the same but the constant $c$ may change.
Similarly, using (ii), $N$ may change while $\beta$ remains the same.
\ \ \qed
\enddemo
 
Thus in particular the $\Cal C^{1+\gamma}$ conjugate of a strictly hyperbolic
set is still hyperbolic, though strictness may be lost.
We remark without proof that a converse holds: by a well-known theorem
due in its original form to Mather, any hyperbolic Cantor set is conjugate 
to a
strictly hyperbolic set (without changing the order of differentiability). 
Therefore if one is studying
properties invariant with respect to the equivalence relation given by 
conjugacy, one might as well 
begin with the assumption that $c=1$ in (ii); this situation
occurs often in the dynamical systems literature. The new metric on $I$ is 
referred to 
as an \it adapted metric\rm~ for the hyperbolic map $S$. 
However for our purposes it will be 
important to use the original metric; otherwise the notion of the
scenery process will lose its meaning.
This will become clear in \S5.
 
We recall from the theory of Bowen and Ruelle (\cite{Bo 1,2},
\cite{Ru}; see also \cite{Be}) that the dimension $d$ of a hyperbolic 
$\Cal C^{1+\gamma}$
Cantor set $C$ is
strictly between 0 and 1 and that {\bf Hausdorff (or conformal)
measure} $\mu = H^d|_C$ has a unique normalized invariant version
$v$, called the {\bf Gibbs measure} (or \it Gibbs state\rm).  For the 
middle-third set
$\mu=v$; in general they are boundedly equivalent (i.e\. the Radon-
Nikodym derivative is bounded away from 0 and $\infty$); $v$ is
defined so as to be a probability measure while $\mu$ may have
total mass $\neq 1$. As in \cite{BF 1}, we will need to use both measures.
 
We remark that all the results in this paper generalize with minor notational 
changes to the following situation:
the maps $\varphi_i$ are also allowed to be orientation-reversing; there 
may be more than two maps,
$\varphi_1,\dots,\varphi_k$, and the Cantor set is constructed by selecting 
the maps with respect
to some subshift of finite type $\Sigma_A$ on $k$ symbols instead of the full 
two-shift $\Sigma.$

\subhead{\S 1.2 Ratio Cantor sets}
\endsubhead
As before, $\Sigma^+$
denotes $\Pi_0^\infty \{0,1\}$
and now we define: $\Sigma^-\equiv \Pi_{-\infty}^{-1} \{0,1\}$,
$\Sigma \equiv \Pi_{-\infty}^{\infty} \{0,1\}$.
We will write
$y=(\dots y_{-2} y_{-1}) $ for
$y\in \Sigma^-$,
$x=(x_0 x_1\dots) \in \Sigma^+$
and
$
\underline x =(y,x)=(\dots x_{-2} x_{-1}.x_0 x_1\dots)=
(\dots y_{-2} y_{-1}.x_0 x_1\dots)
$
for a point in $\Sigma$.
We will let $\sigma$ denote both the  (full) left shift on $\Sigma$
and the left shift (with truncation) on $\Sigma^+$.
$\Sigma^+$
is known as the {\bf future} of $\Sigma$, and
$\Sigma^-$ as its {\bf past}.
 
Write $\Delta$
for the unit simplex in $\r^3$ and int$\Delta$
for its interior.
Let $R$ be a continuous function from $\Sigma^-$ to
int$\Delta$, and write the components $R=(R_l,R_g,R_r)$.
These letters will stand for {\it left, gap} and {\it right}
respectively; by definition they add to $1$ and each is strictly positive.
For $\underline x=(y,x)\in \Sigma$ we will also think of $R$ as a function
on $\Sigma$, by defining $R(\underline x)\equiv R(y)$.
 
Given the function $R$, we will define for each $y\in \Sigma^-$
the {\bf ratio Cantor set}
$C^y\subseteq [0,1]$
so as to satisfy the following:
at each stage, subintervals will have length ratios $R(\sigma^n(\underline x)$.
Thus, we first define
$I^y_0=[0,R_l(y)]$, $I^y_1=[1-R_r(y),1]$.
The left interval $I^y_0$ has subintervals
$I^y_{00}$, $I^y_{01}$
which are defined to have lengths in the ratios
$$
\frac{|I^y_{00}|}{|I^y_{0}|}=R_l (\dots y_{-2} y_{-1} 0.),
\frac{|I^y_{01}|}{|I^y_{0}|}=R_r (\dots y_{-2} y_{-1} 1.)
$$
and left and right endpoints the same as those of $I^y_0$, respectively.
Inductively, for
$x\in \Sigma^+$ and $\underline x=(y,x)$,
$
I^y_{x_0\dots x_{n+1}}$
is a subinterval of
$I^y_{x_0\dots x_{n}}$
with length ratio
$$
\frac{|I^y_{x_0\dots x_{n+1}}|}{|I^y_{x_0\dots x_{n}}|}=
R_*(\sigma^n \underline x).
$$
Here $*=l$, i.e\. this is the left subinterval, if $x_{n+1}=0$, and $*=r$, 
i.e\. the right subinterval, if
$x_{n+1}=1$.
 
Note that the fact that $R$ depends only on the past coordinates 
$y\in \Sigma^-$
of $\underline x$ is what makes this well-defined, since therefore the ratio
is the same for each other point in that subinterval.
 
Now finally we form the set $C^y$ as before, defining
$$
C^y=\bigcap_{n=0}^\infty \bigcup I^y_{x_0\dots x_n}.
$$
 
The simplest example is again the middle-third set: taking 
$R(y)=(\frac 13, \frac 13,
\frac 13)$
for all $y\in \Sigma^-$,
$C^y$ is the middle-third set $C$ for each $y$.

Since $R$ by assumption is a continous map from a compact set into the
 interior
of $\Delta$, hence strictly into the interior, these intervals 
$I_{x_0 x_1\dots}$ nest down to a single
point in $[0,1]$. Hence
each $C^y$ inherits
from $\Sigma^+$ the dynamics of the shift map $\sigma$.
(We will also write $\sigma $ for this map on $C^y$).
 
For $R$ assumed to be H\"older continuous with some exponent $\alpha>0$
(which will be the case in
the present paper), it turns out that each ratio Cantor set $C^y$ is also
a hyperbolic $\Cal C^{1+\gamma}$ Cantor set. (What needs 
to be shown is that $\sigma$
is $\Cal C^{1+\gamma}$ 
on $C^y$ for some $\gamma$ with $0<\gamma<1$, which 
follows by a bounded distortion argument,
and that it
can be  extended to a map $S$ on $I_0^y\cup I_1^y$ without losing any
smoothness and which is hyperbolic. This can be proved from a 
 lemma of Sullivan, (3) in the Introduction above. For
a full proof see \cite{BF 3}). 
Therefore it 
also has a Gibbs measure equivalent to Hausdorff measure. The dimension of
$C^y$ and the Gibbs state, viewed as a measure on $\Sigma^+$,
 are the same for each $y\in \Sigma^-$ for the following reason:
for $y,w\in \Sigma^-$, $C^y$ and $C^w$
have the same
scaling function and hence are $\Cal C^{1+\gamma}$  conjugate.
(See \S 7).

\subhead{\S2 Statement of the problem; Bounded distortion}
 \endsubhead
 
Let C be a hyperbolic Cantor set. We wish to describe the geometry of
the sequence of nested Cantor sets one sees along the
way when zooming down toward a point.

Notation will be as follows: for $x\in C$ with $x=\pi(x_0 x_1\dots)$,
let $C_{n,x}$ be the set $C\cap I_{x_0\dots x_n}$ affinely
rescaled to the unit interval, so as to have endpoints at 0 and 1.
 
Write $\mu_{n,x}$ for the corresponding Hausdorff measure, the restriction
of $H^d$ to the set $C_{n,x}$.  Thus we want
to see how the  sequence of sets, and of the corresponding  measures, varies 
as
$n\longrightarrow\infty$.
 
The first obstacle we encounter is that you don't get from the
interval $I_{x_{0}\dots x_{n}}$ to its subinterval $I_{x_{0}\dots
x_{n+1}}$
by one application of $\varphi_0$ (or $\varphi_1)$.
Instead you have
$$
I_{x_0\dots x_{n+1}}
=\varphi_{x_{0}}\dots
\varphi_{x_{n+1}}(I)=\varphi_{x_{0}}\dots\varphi_{x_{n+1}}
(\varphi_{x_{n}}^{-1}\dots\varphi_{x_{0}}^{-1}(I_{x_0\dots x_{n}})),
$$
and since the maps don't commute, you have to go all
the way back up and down the tree again,
with more nonlinearity introduced each time.
To control this nonlinearity we will use the well-known
Bounded Distortion Property, in the following variation. For a proof
see \cite{Sh-Su}, \cite{Ma\~n\'e} or Lemma 6.4 below. We learned this version 
of bounded
distortion from M. Urbanski.
 
\proclaim{Theorem 2.1 (Classical)}  With $S$ as above,
$\exists K>0$ such that for all $n$,
for any  $\delta >0$,
if $J$ is an
interval such that $S^m|_J$ is 1-1 and the image $S^m(J)$ has
diameter less that $\delta$,
then for all $x,y\in J$,
$$ e^{-K\delta^\gamma}<|\frac{DS^mx}{DS^my}| <
e^{K\delta^\gamma}. $$
\endproclaim
 
We mention that one sees from the proof that if $c$ is the H\"older constant 
for
$\log |DS|$, 
then the constant $K$ is given by
$K=c\beta^\gamma/(1-\beta^\gamma)$.
 
As a consequence of this theorem,
since $0<\alpha <D\varphi_i<\beta<1$ implies that
$\alpha^n<|I_{w_o\dots w_n}|<\beta^n$ for any $w$, we have:
 
\proclaim{Corollary 2.2}  For any $m,n\geq 0$ and any
$w\in\underset o\to{\overset \infty\to \prod}\{0,1\}$ one has
for all $x,y\in I_{w_0\dots w_{n+m}},$
$$
e^{-K\beta^{n\gamma}} < \frac{|DS^m(x)|}{|DS^m(y)|} < e^{K\beta^{n\gamma}}.
$$
\endproclaim

The set $C_{k,x}$ belongs to the collection of $2^k$ Cantor sets at
level $k$ in the tree (rescaled).  We want to understand the geometry of the
sets in this collection.
 
A first approximation is the original set itself (at level 0).  But
by bounded distortion, for $n$ large the $2^n$ sets at level $n$
provide much better models for the $2^{m+n}$ sets at level $k=m+n$;
moreover (and this is the strength of bounded distortion), this
is true for {\it all} $m$ simultaneously.  The reason is that since
by definition
$$
S^m(I_{y_0\dots y_m x_0\dots x_n})= I_{x_0\dots x_n},
$$
and since $I_{x_0\dots x_n}$ has small diameter, by the Corollary the
derivative of $S^m$ is close to constant - hence $S^m$ is close
to linear.

In summary,  consider all the
Cantor subsets which have this same image under $S^m$ to be grouped in one 
equivalence class.
The $2^k$ sets at level $k$ are split into $2^n$
equivalence classes, each with $2^m$ members which all have
approximately the same geometry, (but whose locations are scattered throughout 
the space!). As we scale down toward a point
$x$, we are seeing sets given by these approximations.
 
Note that the
equivalence class of a given interval at level $k$ depends
on the immediately previous $n$ branches, rather than on its
initial branching structure.  We will see in the next section
how Sullivan uses this
observation to study the asymptotics, associating to the Cantor
set a function $R$ like that used to define the ratio Cantor sets in the
previous section. Then, in \S5, we will show that the sequence of sets one 
sees in
$C$ is asymptotically the same as that for the ratio Cantor sets $C^y$.
 
And now for the set $C^y$, the nested sequence of subsets has an exact 
description.
Each subset is itself a ratio Cantor set. Moreover the sequence changes in the 
following way.
Writing as above $C^y_{n,x} $
for the set
$C^y\cap
I^y_{x_0\dots x_n} $ affinely rescaled to $[0,1]$,
one has immediately from the definitions that
$C^y_{n,x} =C^{\sigma^n({\underline x})}$ for all $n\geq 0$.
With the  Gibbs measure on the full shift,
this gives a stationary, set-valued process - which in forward time
describes exactly
what one sees as one zooms down toward Hausdorff-almost every point in
the ratio Cantor set
$C^y$.
 
\subhead{\S3 Sullivan's Scaling Function}
\endsubhead
Now we return to the study of a hyperbolic Cantor set $C$.
Instead of treating the structure of the entire set $I_{x_{0}\dots
x_{n}} \cap C$, which is what we have been emphasising so far, Sullivan 
focuses on the information contained in the 
first step of its
construction, given
by the relative lengths of the subintervals of $I_{x_{0}\dots x_{n}}$.
These subintervals are the left third $I_{x_{0}\dots x_{n}0}$,
right third $I_{x_{0}\dots x_{n}1}$ and middle gap written
$G_{x_{0}\dots x_{n}}$.  We normalize the lengths of these three
intervals, defining for $x\in C $ and
$n\geq0$, where $x=\pi(x_0 x_1\dots)$,
$$R_{n,x}=
(|I_{x_{0}\dots
x_{n}0}|, |G_{x_{0}\dots x_{n}}|, |I_{x_{0}\dots x_{n}1}|)
\big/ |I_{x_{0}\dots x_{n}}|. \quad
$$
This is called by Sullivan the {\bf ratio geometry function} of $C$;
it maps $\Bbb N \times C$ to the interior of the unit simplex
$\Delta\subseteq\Bbb R^3$
and determines $C$ uniquely (one simply constructs $C$ to have these ratios).
 
Next we write, for $y=(\dots y_{-2} y_{-1})$ in $\widetilde{C}\equiv
\underset{-\infty}\to{\overset{-1}\to\prod} \{0,1\}$,
 
$$
R_n(y)=
(|I_{y_{-n}\dots
y_{-1}0}|, |G_{y_{-n}\dots y_{-1}}|, |I_{y_{-n}\dots y_{-1}1}|)
\big/ |I_{y_{-n}\dots y_{-1}}|. \quad
$$

Following Sullivan, it is nice to think of $\widetilde{C}$ as a
distinct Cantor set, dual to $C$ (and called the {\bf dual Cantor set}).
Later for the dynamical interpretation we will instead view
$\widetilde{C}$ as $\Sigma^-$,
 that is,
as the past coordinates of the full shift $\Sigma
\equiv \underset{-\infty}\to{\overset \infty\to\prod} \{0,
1\}$. We will use whichever symbol ($\widetilde{C}$ or $\Sigma^-$)
is more appropriate in the context.


As in \cite{Bo 1}, for any $\beta\in (0,1)$, the $\beta$-metric on $\Sigma^-$
 (which defines what is
meant below by H\"older continuity) is taken to be:
$d_\beta(y,w)=\beta^n$ where $n$ is the greatest
positive integer such that $y_{-n}\dots y_{-1} = w_{-n}
\dots w_{-1}$. When $\Sigma^-$ is thought of as $\widetilde C$, i.e\. as dual 
to a specific hyperbolic Cantor set $C$,
we choose $\beta$ to be (as before) the upper bound on $D\varphi_i$.
We mention that if $\beta$ is replaced by some other number
$\tilde\beta\in (0,1),$ then the metrics are related by
$$d_{\tilde\beta} =(d_\beta)^{\log\tilde\beta/\log\beta},
$$
and the H\"older exponent $\gamma$ for $R$ in the statement of the next 
theorem would change
to $\gamma\cdot(\log\tilde\beta/\log\beta)$.

\proclaim{Theorem 3.1 (Sullivan)}  Let $C$ be a hyperbolic 
$\Cal C^{1+\gamma}$ Cantor
set. For every $y$ in the dual Cantor set $\widetilde C$,
$$
R(y)\equiv\lim_{n\to\infty}
R_n(y)
$$
exists. The convergence is of order $O(\beta^{n\gamma})$, uniformly in $y$, 
and
the function $R$ is H\"older
continuous with exponent $\gamma$, in the $\beta$-metric. $R$ takes values 
strictly in the interior of $\Delta$. 
\endproclaim

\proclaim{Definition} 
$R$ is called the {\bf scaling function} of $C$.
\endproclaim 
 
\demo{Proof}  We will first show that for each $y,
$ $R_n(y)\ \ n=1,2,\dots$ is a Cauchy sequence.  Since
$$
S^m(I_{y_{-(n+m)}\dots y_{-1}})=I_{y_{-n}\dots y_{-1}}
$$
and similarly for the subintervals, applying the Mean Value Theorem
and Bounded Distortion Property (Corollary 2.2) we have for all $m\geq 0$
$$
R_n(y)=R_{n+m}(y) e^{\pm K\beta^{n\gamma}}
$$
Therefore $R_n(y)$ is  Cauchy sequence (i.e\. each of its three
coordinates is) hence it converges; call the limit
$R(y)$.  Next, if $y,w\in
\underset{-\infty}\to{\overset{-1}\to\prod}\{0,1\}$ agree on the
coordinates $-n,\dots,-1$ then since
$R_n(y)=R(y)e^{\pm K\beta^{n\gamma}}$ and
$R_n(y)=R_n(w)$, we have
$$
R(y)=R(w) e^{\pm 2K\beta^{n\gamma}}.
$$
Writing $\parallel\cdot\parallel$ for sup norm in ${\Bbb R}^3$, this
implies that, with the log taken by components,
$$
\parallel\log R(y)-
\log R(w)\parallel\leq
2K(d_\beta(y,w))^\gamma,
$$
i.e\. $\log R$ is H\"older continuous with exponent $\gamma$;
therefore so is $R$.
\ \ \qed \enddemo

Thus (since the normalized lengths add to one), $R$ maps $\widetilde C$
onto a compact subset of the interior of the
unit simplex in ${\Bbb R}^3$.
 
\subheading{\S 4 Dynamical versions of Sullivan's theorem}
In this and the next section we return to the original motivating
question: what
does the sequence of sets $C_{n,x}$ look like?  This
is exactly what one sees for the $n^{\text{th}}$ level Cantor set, as one
zooms down toward $x$.
 
First we state Theorem 3.1 in a dynamical form. Here it will be crucial to
think
of the dual Cantor set $\widetilde C$ as the past $\Sigma^-$ of $\Sigma$.
We extend the function $R$ to $\Sigma$ by defining:
$R(\underline x)=R(y)$ for $\underline x= (y,x)$.
This function depends only on the past coordinates y of $\underline x$.
 
\proclaim{Corollary 4.1}  For each $x\in C$, for any choice of
$\underline{w}\in\Sigma$ such that $\exists k \geq 0$ with
$w_k,w_{k+1}\dots =x_k,x_{k+1}\dots$ (where $x=\pi(x_0 x_1\dots)$),  then
$$\parallel R_{n,x}-R(\sigma^n\underline{w})\parallel
\longrightarrow 0 \quad\text{as $n\longrightarrow\infty$}.
$$
\endproclaim
 
The proof is immediate from the definitions, and in fact if the
H\"older constant for $R$ is $c>0$ so that
$$\parallel R(\underline{x})-R(\underline{z})\parallel \leq c
d_\beta(\underline{x},\underline{z}),
$$
one has
$$
\parallel R_{n,x}-R(\sigma^n(\underline{w}))\parallel\leq
c \beta^{n-k};
$$
here the $\beta$-metric has been extended to $\Sigma$ in the natural
way, with points $\underline x,$
$\underline w$ having to agree on coordinates from $-n$ to $n$.
 
We note that equivalently,
if $x$ and $w$ are in the same unstable set in $(C,S)$ then
the sequence $R_{n,x}$ (with any past)
is in the stable set (in the shift on sequence space)
of the sequence given by $R$ sampled along the shift orbit of
$\underline x$, with
an exponential rate of convergence.

We recall that a \it stochastic process \rm is simply a (one- or two-)
sided sequence of measurable functions $f_i$ (known as \it
random variables\rm) defined on some probability space $(\Omega, \nu)$. 
The process is
\it stationary \rm if a time-shift doesn't alter the probability of an event.
Equivalently, the space of paths 
$\{(\dots,f_i(\omega),\dots)\}$, acted on by the shift tranformation and given
the pushed-forward measure, is a measure-preserving transformation of a 
probability space. Conversely, a 
measure-preserving transformation determines many 
stochastic  processes: choose a measurable function and evaluate it along
orbits. Thus for example $R(\sigma^n\underline{w})$ is a
($\Delta$-valued) stochastic process. In the next section we will
encounter set- and measure-valued versions of this.

Next, recall the definition of a {\bf generic point}  $x$ for an ergodic
measure-preserving transformation $T$ on a compact metric space $X$ with
probability measure $m$.
For each continuous $f:\ X \to \r$, $x$ satisfies:
$$
\lim_{N\to\infty} \frac 1N \sum_{k=0}^{N-1} f(T^kx)=\int_X
f\ dm.
$$
That is, $x$  samples each continuous function well with respect
to time averages. If $X$ is a Polish space (a complete
separable metric space) --
this will occur in the next section -- ,
then we instead sample the continuous functions with compact support.
By the remarks in the previous paragraph,   this definition also makes
sense for a stationary ergodic stochastic process, 
 if the path space has been given the topology of  a Polish space.
In the definition of generic point we only take time averages toward 
$+\infty$, so in the case of a two-sided stochastic process, 
it will be natural to allow a
  one-sided sequence as a  generic point, as well.

If the measure lives in a  compact part of the space (which will always be 
the case in this paper) then
 by an ergodic theorem of Kryloff and Bogliouboff (i.e\. by the Birkhoff ergodic
theorem plus compactness), $m$-almost every $x$ is a  generic point.
 
Now $\mu$ is equivalent to the Gibbs measure $\nu$, which is
invariant and has a unique invariant natural extension $\hat{\nu}$
on $(\Sigma, \sigma)$.  Hence by Kryloff
and Bogliouboff:

\proclaim{Corollary 4.2}  For $\mu-$ a.e\. $x\in C$, the (one-sided) sequence
$R_{1,x}, R_{2,x}\dots$ is a generic point for the ergodic
$\Delta$-valued process $R(\sigma^n(\underline{w}))$ for
$n\in{\Bbb Z}$, given by $\underline{w}\in\Sigma$ being distributed like
$\hat{\nu}$.
\endproclaim

\subheading{\S 5 Conjugacies, and the scenery process}
 In this section we will construct a set-valued version of the 
scaling function, and use it  to prove analogues of Corollaries 4.1 and 4.2,
which will describe how 
 the sequence of sets $C_{n,x}$ approximates
 the scenery process.
We will make use of three different metrics on collections
 of Cantor sets.  One metric, which is
derived from the 
 $C^1$ norm on the space of \it conjugacies \rm of Cantor sets, 
is well suited to proofs and  is natural 
from an abstract point of view.  There  we will prove properties
(convergence at at exponential rate; H\"older dependence) which
will then  pass over to  two geometrically defined metrics: the 
Hausdorff 
metric, and a metric derived from the Hausdorff measures.

\subheading{Three metrics}
Our  metrics will be defined on several different spaces.
Fix a hyperbolic $\Cal C^{1+\gamma}$
Cantor set  with map, $(C,S)$.
We write $\Cal E^{1+\gamma}
\equiv \Cal E^{1+\gamma}(C)$ 
for the collection of 
Cantor sets (with maps)
which are $\Cal C^{1+\gamma}$-conjugate to $(C,S)$. (From Lemma 1.1.2, 
these are also 
hyperbolic $\Cal C^{1+\gamma}$ Cantor sets). 
We write 
$\Cal E^{1+\gamma}_*$ 
for the quotient space of $\Cal E^{1+\gamma}$ 
where $(C,S)$ and $(C,\wt S)$ 
are identified if  $S=\wt S$ on $C$.
This is the collection of marked Cantor sets conjugate to $C$, or 
equivalently the pairs $(C,S|_C)$ with restricted maps.
$\Cal E^{1+\gamma}_{**}$ 
will denote the collection of underlying Cantor sets.
We write $\text{Diff}^{1+\gamma}$
for the $\Cal C^{1+\gamma}$ order-preserving diffeomorphisms of $I$.
Given choice of the pair $(C,S)$, $\text{Diff}^{1+\gamma}$
projects onto  $\Cal E^{1+\gamma}(C)$ 
in a natural way:
$f$ is mapped to $(C_f,S_f)\equiv     (f(C), f\circ  S\circ f^{-1})$.
This is many-to one because there is some freedom given by the gaps;
see Proposition 8.3. We note that  the projection from 
$\Cal E^{1+\gamma}$
to $\Cal E^{1+\gamma}_*$ 
is   also many-to-one.

The scenery process can be thought of as taking values 
in $\Cal E^{1+\gamma}_{*}$, 
the marked Cantor sets, or in the space of 
underlying sets $\Cal E^{1+\gamma}_{**}$.
We will first prove convergence in the space of conjugacies, 
$\text{Diff}^{1+\gamma}$; this will then imply convergence in
 the other spaces.

First we consider two metrics on $\Cal E^{1+\gamma}_{**}$.
We recall the definition of  the {\bf Hausdorff metric} on the collection
of  closed subsets
 of the interval $I$:
$$d_H(A,B)=\text{inf}  \{\epsilon:A+(-\epsilon,\epsilon)\supseteq
B\text{ and }B+(-\epsilon,\epsilon)\supseteq A\}.$$
This defines a metric on  $\Cal E^{1+\gamma}_{**}$, and a 
pseudo-metric on the other spaces defined above.

Next, we define the following metric  
on 
 the set of finite Borel measures on [0, 1], denoted $\Cal M$.
Enumerating binary  intervals $E_1, E_2,
\dots, E_n\dots$ of the form $[j2^{-k},(j+1)2^{-k}]$,  for $\nu_1,\ \nu_2$
in $\Cal M$, set
$$d(\nu_1,\ \nu_2)=
\sum_{n=1}^\infty |\nu_1(E_n)-\nu_2(E_n)|/2^{-n}.$$
 This metric induces a topology equivalent to 
the weak topology on $\Cal M$, in the language of  probability theory;
in  analysis terminology this is  
the weak-$*$ topology on $\Cal M$, the dual of 
the space of continuous functions. 

On $\Cal E^{1+\gamma}_{**}$ we define the {\bf measure metric} $d_M$
from this,  setting: 
$$d_M(C,D)=
d(H^d |_C ,H^d |_D)$$ where $H^d$ is  $d$-dimensional Hausdorff measure.
On the other spaces, this again defines a  pseudo-metric.

Next, recall that the $\Cal C^1$-norm
of $f:I\to \r$ is:
$$
\|f\|_{\Cal C1}=\|f\|_{\infty}+\|Df\|_{\infty}.
$$
 We identify  $\text{Diff}^{1+\gamma}$ 
with the collection of triples
$(C_f, S_f, f)$ for $f\in \text{Diff}^{1+\gamma}$, to be written as
$\wh {\Cal E}^{1+\gamma}$.
(As we noted above, the map from
$\wh {\Cal E}^{1+\gamma}$
to $\Cal E^{1+\gamma}$ 
is not one-to-one). 
The $\Cal C^1$-norm on $\text{Diff}^{1+\gamma}$ determines a metric
on $\wh {\Cal E}^{1+\gamma}$ as follows.
For $f,\ g$ in $\text{Diff}^{1+\gamma}$, we write 
$$d_C(C_f,C_g)= \|f-g\|_{\Cal C1}. $$
 
We  call this the {\bf $\Cal C^1$ metric} on 
$\wh {\Cal E}^{1+\gamma}$. 

\subhead{Note}
\endsubhead
The metric $d_C$ keeps track of the map $S$ on all of its domain
$I_0\cup I_1$, 
while $d_H$ and $d_M$ only see the Cantor sets.
Thus for the $\Cal C^1$ metric, writing $\Cal I$
for the identity map,
 $(C,S, \Cal I)$ and $( C,S_f,f)$ 
 will be a positive distance apart unless  in particular
$S= S_f$ on all of $I_0\cup
 I_1$. (Of course also one needs $f=\Cal I$).  
\smallskip

The definition of $d_C$ depends on the initial choice of the set $(C,S)$ (with
the identity map);
the next proposition shows how the metric varies if we change this  ``base 
point'' of $\wh {\Cal E}^{1+\gamma}$. 
As we will see in   \S8, $\text{Diff}^{1+\gamma}$ is a Lie group
and this is also a  statement about bounded invariance of 
a metric on that group.

\proclaim{Proposition 5.1}  Let $D\in \Cal E^{1+\gamma}$, with 
$D=\Phi(C)$, with $\Phi\in \text{Diff}^{\ 1+\gamma}$. We have:

$$\frac 1K d_D< d_C < K d_D$$
 where $K=2\text{max} \{ \|D\Phi\|_\infty, \|D(\Phi^{-1})\|_\infty^{-1} \}$.
\endproclaim
\demo{Proof}
We note that from the definition of the  $\Cal C^1$-norm
one has that if $f:I\to I$ with $f(0)=0$, then
$
\|f\|_{\Cal C1}\leq 2\|Df\|_{\infty}.
$
Therefore if also $g(0)=0$, then
$\|g\circ f\|_{\Cal C1}\leq 2\|g\|_{\Cal C1} \|f\|_{\Cal C1}.$

Now for $f,\ g$ and $\Phi$ as in the statement of the Proposition,
we have:
$$
\aligned
d_D(C_f,C_g)
&=\|f\circ\Phi^{-1}-g\circ\Phi^{-1}\|_{\Cal C1}
=\|(f-g)\circ\Phi^{-1}\|_{\Cal C1} \\
&\leq 2\|f-g\|_{\Cal C1}\|\Phi^{-1}\|_{\Cal C1}
=d_C(C_f,C_g)
\|\Phi^{-1}\|_{\Cal C1}
\endaligned
$$
which gives one of the inequalities. The other is proved in the same way.
\ \ \qed
\enddemo

Next we will look at how (on $\wh {\Cal E}^{1+\gamma}$) 
the pseudo-metrics $d_H$ and $d_M$ compare to the metric $d_C$.
First we recall how
the Hausdorff measure transforms under mappings.

\proclaim{Definitions}
Given a 1-1 differentiable map
$\varPsi:M\longrightarrow N$ between open subsets of ${\Bbb R}$ and
given a Borel measure $\mu$ on $M$ and real number $d>0$ we write:
$(\varPsi^\star\mu)(E)=\mu(\varPsi^{-1}E)$ and
$$(\check\varPsi \mu)(E)=\int_{\varPsi^{-1}E} |D\varPsi|^d d\mu.$$
Thus $\varPsi^\star\mu$ is the usual {\bf push forward} of $\mu$,
and $\check\varPsi \mu $ is the $(\varPsi,d)-$ {\bf conformal transform}
of $\mu$.
\endproclaim

Hausdorff measure has the {\bf conformal transformation
property} with respect to $\Cal C^1$ 
maps:
for $\varPsi:\r\to \r$ a $\Cal C^1$ 
diffeomorphism,
$$H^d=\check\varPsi(H^d).$$
\proclaim{Proposition 5.2} 
With $d_H$, $d_M$ and $d_C$ denoting the Hausdorff, measure, and $\Cal C^1$ 
metrics respectively, 
for all $C_f,\ C_g \in \wh {\Cal E}^{1+\gamma}$,
$$
d_H(C_f,C_g)
\leq  d_C(C_f,C_g)$$
and for $\Psi(x)=5x+4x^2$, we have for all $C_f\in \Cal E^r(C)$,
$$
d_M(C,C_f)
\leq \Psi(d_C(C,C_f)).$$
\endproclaim
\demo{Proof}
For the Hausdorff (pseudo)-metric this is
immediate, using
the $ L^\infty$ norm, since
$$d_H(C_f,C_g)\leq \|f-g\|_\infty\leq
 d_C(C_f,C_g).$$

For the second inequality, writing $\Cal I$ for the identity map on
$I$, we have
$$d_C(C,C_f)\equiv\|f-\Cal I\|
_{\Cal C 1}$$
and writing $\mu=H^d|_C, \ \mu_f=H^d|_{C_f},$
$$
\aligned
d_M(C,C_f)
&\equiv\sum_n|\mu E_n-\mu_f E_n|  2^{-n}\\
&=\sum_n\bigl |\int_{f^{-1}E_n}|Df|^d\ \text{d}\mu-\int_{E_n}1\ \text{d}\mu
\bigr | 
2^{-n}\\\\
&\leq 2 \|f-\Cal I\|_\infty 
\|Df\|_\infty^d 
\sum_n 2^{-n}  
+ \sum_n 2^{-n}\int_{f^{-1}E_n
\cap E_n}
\bigl | |Df|^d-1\bigr |
\ \text{d}\mu
\endaligned
$$
(here the first term bounds the contribution, for each interval $E_n$,
of its two ends not matching up exactly with $f^{-1}E_n$; we used the fact that
since $f$ is a diffeomorphism of $I$, $\|Df\|_\infty\geq 1$).

Next, we note that for all $x>0$, $|x^d-1|\leq |x-1|.$
Hence $\bigl||Df|^d-1\bigr|\leq \bigl||Df|-1\bigr|,$ so the above is
$$\aligned
&\leq
2\|f-\Cal I\|_\infty(1+\bigl\||Df|-1\bigr\|_\infty) +
\bigl\||Df|-1\bigr\|_\infty \sum_n |E_n|
2^{-n} 
\\
&\leq
4\bigl\||Df|-1\bigr\|_\infty (1+\bigl\||Df|-1\bigr\|_\infty) +
\bigl\||Df|-1\bigr\|_\infty 
\\
&\leq
4(1+\bigl\||Df|-1\bigr\|_\infty+1)( \bigl\||Df|-1\bigr\|_\infty)
=5\bigl\||Df|-1\bigr\|_\infty+ 4(\bigl\||Df|-1\bigr\|_\infty)^2 
\\
&\leq\Psi( d_C(C,C_f)),
\endaligned$$
as claimed.
\ \ \qed\enddemo

We define the {\bf (restricted) $\Cal C^1$  metric} on $\Cal E^{1+\gamma}_*$ 
to be 
$$d_{C*}
(C_f,C_g)= \|f-g\|_{\Cal C1, C}$$
where this indicates that the sup norms are taken over the Cantor set $C$. 
Obviously $d_{C*}\leq d_{C}$ on $\wh{\Cal E}^{1+\gamma}$.

The first theorem we are aiming for will state:

\proclaim{Theorem 5.3}  Given a $\Cal C^{1+\gamma}$ hyperbolic Cantor set 
$(C,S)$,
for every $y=(\dots y_{-n}\dots y_{-1})\ \in 
\Sigma^-$,
the limit
$$
C(y)\equiv \lim_{n\longrightarrow \infty} C_{y_{-n}\dots
y_{-1}}
$$
exists. Convergence is exponentially fast, 
and $C(y)$ 
is H\"older continuous with exponent $\gamma$
(for both metrics on $\Cal E^{1+\gamma}_{**}$, 
and for the  restricted $\Cal C^1$  metric on the marked Cantor sets 
$\Cal E^{1+\gamma}_*$). 
Moreover, $C(y)=C^y$, the ratio Cantor set built 
from the scaling function of $(C,S)$.
\endproclaim

Thus $C(y)$ is a set-valued version of the scaling function $R(y)$.
To prove Theorem 5.3, we first construct certain conjugacies (Theorem 5.9),
proving convergence in the $\Cal C^1$ metric.
Convergence
in the Hausdorff and  measure metrics  then follows from 
Proposition 5.2.
Before giving this construction, we state several further consequences
of Theorem 5.9.

In the same way as for the scaling function $R(y)$, 
we define $C(\cdot)$ also on the full shift space $\Sigma$, 
setting
$C(\underline{x})=C(y)$ for $\underline{x}=(y,x)$. We will have:

\proclaim{Theorem 5.4}  For every $x\in C$, the sequence
$C_{n,x}\equiv C_{x_0\dots x_n}$ 
is asymptotic to $C ({\sigma }^n 
(\underline{x}))$,
with an exponential rate of 
convergence, for any $\underline{x}$ in $\Sigma$ with the same future coordinates
$(x_0,x_1,\dots)$.
\endproclaim

We will write $\Cal L_C$ for the (compact) subset of $\Cal E^{1+\gamma}_{*}$
(with respect to all three metrics)
which is the range of the function $\underline{x}\mapsto C({\underline{x}})$. This is the collection of {\bf (marked) limit sets}.
Since the function is continuous and the domain $\Sigma$ is compact,
we then have:

\proclaim{Proposition 5.5} The collection of limit sets is compact.
\endproclaim

\proclaim{Definition} Given a $\Cal C^{1+\gamma}$ hyperbolic Cantor set $C$,
the {\bf set-valued scenery process} of $C$ is the process
$C(\sigma^n (\underline x))=C^{\sigma^n (\underline x)}$, 
with $\underline x\in\Sigma$, distributed according to the measure $\hat\nu$.
\endproclaim

Note that stationarity and ergodicity of this process
follow immediately from invariance and ergodicity of the measure $\hat\nu$.

The space of paths is a compact subset of the Polish
space $\Pi_{-\infty}^{+\infty}
\Cal E^{1}_{**}$ 
(or $\Pi_{-\infty}^{+\infty} \Cal E^{1}_{*}$ for the marked sets),
with the product topology determined by the topologies of any of the 
three metrics. See the proof of Corollary 5.6.

We mention why we use probability terminology -- the \it scenery process \rm--
for the map $C(\underline{x})\mapsto C(\sigma \underline{x})$.  
Note that this  dynamics is not in fact given by a map on 
$\Cal L_C$ itself. Indeed  
at every stage you have two choices -- the right- or left-hand
subsets from the next level 
of the ratio Cantor set $C(\underline{x})$,
with the choice of left or right depending on whether $x_0$ is $0$ or $1$.
Or, from a different viewpoint, one has
 the dynamics  on $\Cal L_C$  of a semigroup 
action; see  the note at the end of the paper.

The next result is like Corollary 4.2. 

\proclaim{Corollary 5.6}  For $\mu-$a.e.$\ x\in C$, the sequence
of (rescaled) Cantor sets $C_{n,x}$ which nest to $x$ 
is a generic point for the stationary ergodic
set-valued process
$C^{\sigma^n (\underline x)}$ 
determined by $(\Sigma, \sigma, \hat{\nu})$.
\endproclaim

In \S2, we defined $\mu_{n,x}$ to be the sequence of Hausdorff measures
$H^d|_{C_{n,x}}$. Convergence of $C_{n,x}$ to $C(y)$
in the measure metric can be rephrased as follows:

\proclaim{Corollary 5.7}  For every $y=(\dots y_{-n}\dots y_{-1})\ \in 
\Sigma^-$,
the limit
$$
M(y)\equiv \lim_{n\longrightarrow \infty} \mu_{y_{-n}\dots
y_{-1}}
$$
exists, and is H\"older continuous with exponent $\gamma$.
\endproclaim
 
The support of the measure $M(y)$ is the set $C(y)$, and
$M(y)$ is a measure-valued version of the scaling function.

As we did for the set-valued process, we 
define $M(\underline{x})=M(y)$.
We define 
the
{\bf measure-valued scenery process}
$M(n,\underline x)\equiv M(\sigma^n (\underline x))$,
again with
$\underline x\in\Sigma$, distributed like $\hat\nu$.
We have:

\proclaim{Corollary 5.8}  For $\mu-$a.e.$\ x\in C$, the sequence
$\mu_{x_0\dots x_n}$ is a generic point for the stationary ergodic
measure-valued process
$
M( \sigma^n(\underline{x}))$
determined by $(\Sigma, \sigma, \hat{\nu})$.
\endproclaim

Now we proceed to the proofs.
For each interval $I_{ w_0\dots w_n}$ with $ w_i=0$
or 1, we write $A_{ w_0\dots w_n}$ for the affine map which
expands the interval to the unit interval I.  
We then set for $y\in \wt C$,
$A^y_n= A_{y_{-n}\dots y_{-1}}.$
This expands the interval
$I_{y_{-n}\dots y_{-1}}$
to $I$ affinely. Next, define maps
$\varphi^y_{k,n}:I\to I$
for $k\geq n$ by:
$$
\aligned
&\varphi^y_{k,k}= \text{ identity and}\\
&\varphi^y_{k,n}= \varphi_{y_k}\circ\dots\circ\varphi_{y-(n+1)}\\
\text{ for } k\geq n.
\endaligned
$$
We will also write $\varphi^y_k$
for
$\varphi^y_{k,0}.$

For $0\leq n\leq k$, we define
$ \Phi^y_{k,n}:I\to I$
by:
$$
\Phi^y_{k,n}= A^y_k
\circ\varphi_{k,n}^y\circ (A^y_n)^{-1}.
$$
It follows that:
$$
\align
\text{ for each }  n,\phantom{elephant}
&\Phi^y_{n,n}=
\text{ identity,} \tag1 \\
\text{ for all $m\geq k\geq n$,} \tag2 \\
          &\Phi^y_{m,n}= \Phi^y_{m,k}\circ\Phi^y_{k,n}
          \text{ and }\\
          &\Phi^y_{n,0}= A^y_n\circ \varphi^y_n. \tag3
\endalign
$$
We will also write $\Phi^y_n$
for $\Phi^y_{n,0}.$
 
The sequence
$\Phi^y_n$
is, to use Sullivan's words, ``a sequence of ... compositions (of contractions)
... renormalized by post composition with linear maps
to obtain mappings between unit intervals...''
(Appendix of \cite{Su 1}). As Sullivan states, and as we will prove in the next
section, such a sequence
is precompact in $\Cal C(k,\alpha)$
if the original hyperbolic Cantor set is $\Cal C(k,\alpha)$.
 
This gives convergence along some subsequence $\Phi^y_{n_k}$.
However unfortunately, that is not enough for our goal of
proving an ergodic theorem for Cantor sets  and
measures (i.e\. Corollaries 5.6, 5.8) --
for that purpose we want instead to prove that the sequence itself converges.
We do this in the next theorem, using bounded distortion, and
then in the next section we return
to Sullivan's idea to prove smoothness of the resulting limiting conjugacy.
 
\proclaim{Theorem 5.9} Let $(C,S)$ be a
$\Cal C^{1+\gamma}$-hyperbolic Cantor set.
For each $y$ in the dual Cantor set $\wt C$,
$$
\Phi^y\equiv\lim_{n\to\infty} \Phi^y_n
$$
exists. This is an order-preserving diffeomorphism from $I$ to $I$.
Convergence is of order $O(\beta^{n\gamma}) $ in the $\Cal C^1$ norm,
uniformly in $y$, and the function $y\mapsto \Phi^y$ is
H\"older continuous of order $\gamma$, in the $\beta$-metric.
\endproclaim
 
\demo{Proof}
We will show, using bounded distortion, that for $n$ large and for $k>n$
arbitrary, $\Phi^y_{k,n}$
is close to the identity. Then since $\Phi^y_{k,0}=
\Phi^y_{k,n}\circ \Phi^y_{n,0}$,
this will imply convergence.
 
Now since the maps $A$ are affine they have constant derivative.
So for each $a\in I$, for $m=k-n$, we have:
$$
\align
(D(\Phi^y_{k,n}))(a)
\equiv
D(A^y_k\circ (S^m)^{-1}\circ (A^y_n)^{-1})(a)
& =\frac {DA^y_k}{DA^y_n DS^m(z)} \text{ for } z\equiv (A^y_n)^{-1}(a). \\
& =\frac {DS^m(z_0)}{DS^m(z)} \text{ for some }z_0\in I_{y_{-k}\dots y_{-1}},
\endalign
$$
by the Mean Value Theorem.
Therefore by Corollary 2.2
$$
e^{-K\beta^{n\gamma}}
<D\Phi^y_{k,n}(a)
<e^{K\beta^{n\gamma}}
$$
for all $k>n$, $y\in\wt C$ and all $a\in I$.
This implies the sequence $D\Phi^y_{n}
,\ n=0,1,2\dots$
is Cauchy, hence converges. Since
$\Phi^y_{n}(0)=0$
for all $n$,
it follows from the Fundamental Theorem of Calculus that the limit
$\Phi^y\equiv\lim_{n\to\infty}\Phi^y_{n}$
exists, and that
$D\Phi^y=\lim D\Phi^y_{n}$.
By Corollary 2.2,
$D\Phi^y$
is bounded away from $0$ and $\infty$
by
$e^{\pm K\beta^{\gamma}}$;
in particular,
$\Phi^y$
is an order-preserving diffeomorphism from $I$ to $I$, as claimed.

We define for each n
$$
\Phi^y_{\infty,n}=\lim_{k\to\infty}\Phi^y_{k,n};
$$
the limit exists by the above arguments, and this map is
$\Cal C^1$-close to the identity map $\Cal I:I\to I$.
We have for each $n$ that
$$
\Phi^y=\Phi^y_{\infty,n}\circ \Phi^y_{n,0},
$$
and that:
\item{(1)} $\exists k_0>0$ such that
$$\| \Phi^y_{\infty,n}-\Cal I\|_{\Cal C1}<k_0\beta^{n\gamma}.
$$

The constant  $k_0$ here only depends on $K$ from Corollary 2.2,
which in turn depends on $\beta$, the upper bound for
$|D\varphi_i|$. Here is the calculation:
we have
$$\| \Phi^y_{\infty,n}-\Cal I\|_{\Cal C1}\leq
\| D\Phi^y_{\infty,n}-1\|_{\infty}$$
and we know that
$$
e^{-K}\leq e^{-K\beta^{n\gamma}}
\leq D\Phi^y_{\infty,n}
\leq e^{K\beta^{n\gamma}}\leq e^K
$$
for all $n$.
Now since for $x$ in the interval $[e^{-K},e^{K}]$ one has
$|e^x-1|<k_0x +1$, where
we take
$k_0=(\text{exp}(\text{exp}(K))-1)/(\text{exp}(K))$,
statement $(1)$ holds true.

Now recall from the proof of Proposition 5.1 that
for
 $f,\ g:I\to I$ with $f(0)=0$ and $g(0)=0$ then
$\|g\circ f\|_{\Cal C1}\leq 2\|g\|_{\Cal C1} \|f\|_{\Cal C1}.$
From this, it follows that statement $(1)$ is equivalent to:

\item{(2)} $\exists k_1>0$ such that
$$\| \Phi^y_{n}-\Phi^y\|_{\Cal C1}<k_1\beta^{n\gamma}.
$$
(Here we can take $k_1=(e^K)k_0).$

Both statements express, in different ways, that $\Phi^y $ is close to
$\Phi^y_{n,0}$, with exponentially fast
 convergence; $(2)$   is what we stated in the Theorem.
Finally it is now also easy to check H\"older continuity:
$\exists
 k_2>0$ such that
for all $y,w\in\wt C$,
$$
\| \Phi^y-\Phi^w\|_{\Cal C1}\leq k_2 (d_{\beta}(y,w))^{\gamma}.
$$
\phantom{00000000000000000000000000000000000000000000000000}\ \ \qed \enddemo

\demo{Proof of Theorem 5.3}
Writing $C^y_n$
for $\Phi^y_n(C)=C_{y_{-n}\dots  y_{-1}}$,
and $C^y\equiv\Phi^y(C)$,
statement $(2)$ above says exactly:
$$d_C(C^y_{n},C^y)<k_1\beta^{n\gamma}.$$ Hence we have convergence in the 
$\Cal C^1$ metric.
By Proposition 5.2 therefore, $d_H$ has  the same bound. 
For the measure metric, we have 
$$\aligned
d_M(C^y_{n},C^y)
&\leq \Psi(d_{C^y}(C^y,C^y_{n}))\\
&=
\Psi(\| \Phi^y_{\infty,n}-\Cal I\|_{\Cal C1})\\
&\leq\Psi(k_0\beta^{n\gamma}) 
\endaligned
$$
by $(1)$, and this is
$\leq k_3\beta^{n\gamma} 
$ where $k_3=k_0(5+4k_0).$
Next we show H\"older continuity.
Now $d_C(C^y,C^w)
=\| \Phi^y-\Phi^w\|_{\Cal C1}$
so $(3)$ proves H\"older continuity for the $\Cal C^1$ and Hausdorff metrics.
Then, applying Proposition 5.1,
$d_M(C^y,C^w)\leq \Psi(d_{C^y}(C^y,C^w))
\leq\Psi(d_C(C^y,C^w)\|\Phi^{-1}\|_{\Cal C1})
\leq \Psi(e^Kk_2 \beta^{n\gamma})\leq k_4 \beta^{n\gamma}$
where 
$k_4=5a+4a^2$ and $a=e^Kk_2.$
Finally, it is clear from the constructions that $C(y)$ 
has ratio geometry 
given by $R(y)$, hence $C(y)$  is indeed equal to $C^y$.
\ \ \qed\enddemo

\demo{Proof of Theorem 5.4} 
From the proof of Theorem 5.3,
since the exponential bound is uniform over all sets of level $n$,
we have that the $d_C$-distance from 
$C(\sigma^n x)$ to $ C_{x_0\dots x_n}$ is bounded by 
$k_1\beta^{n\gamma}$.
The bounds for $d_H$ and $d_M$ then follow as above.
\ \ \qed\enddemo

\demo{Proof of Corollary 5.6}
We give the proof for the  process which takes  values in the collection of 
marked sets.
Here the space will be $\Pi_{-\infty}^\infty \Cal E^1_{*}
(C)$,
with the shift transformation; this is a Polish space
(since $\Cal E^1_{*}$ is;  we use 
$\Cal E^1_{*}$ rather than $\Cal E^{1+\gamma}_{*}$ 
so as to have a complete space).
Therefore we know from our definitions what it will mean for the one-sided 
sequence 
$C_{n,x}$ 
to be a generic point. Now the map
$\underline x 
\mapsto 
(\dots C(\sigma^{-1}\underline x), C (\underline x),
C(\sigma \underline x)\dots)$
from $\Sigma $
has as its image a compact invariant subset of
$\Pi_{-\infty}^\infty \Cal L_C\subseteq \Pi_{-\infty}^\infty \Cal E^1_{*}(C)$;
this image is the space of paths of the scenery process, and is the support of
the image of the measure $\hat \nu$.
The ergodicity of
$(\Sigma,\hat \nu,\sigma)$  passes over to the scenery process, hence
a.e\. path
$(\dots C(\sigma^{-1}\underline x), C (\underline x),
C(\sigma \underline x)\dots)$
is generic for the shift on path space.
Finally, since by Theorem 5.4 we know the sequence $C_{n,x}$ 
is forward asymptotic to $C(\sigma^n \underline x)$,
we will compute the same time average for the continuous functions.
Thus $C_{n,x}$ is generic, for $\nu$-a.e\. $x$ and hence for 
$\mu$-a.e\. $x$.
\ \ \qed\enddemo

\subheading{\S 6. Smoothness of conjugacies}
Now we will see how to prove the conjugacies of \S 5
in fact have higher smoothness properties. 
 
The basic idea will be to imitate what one knows about analytic maps, for
$\Cal C^{k+\gamma} $
or $\Cal C(k,\gamma)$
maps.
Thus Lemma 6.2 is a version of Leibnitz' formula, and Lemma 6.3 is
one step in showing $\Cal C^{k+\gamma} $ maps are morphisms in a category.
This means they can be used to define equivalence relations on sets,
and to give the analogue  of differentiable structures.
One also imitates the Arzela-Ascoli Theorem, in Lemma 6.4;
as Sullivan says in the Appendix of \cite{Su 1}, and quoted in the previous 
section,
the sequence $\Phi^y_{n}$ will be
precompact in $\Cal C(k,\alpha)$.
(As usual, we do our proofs instead in $\Cal C^{k+\gamma}$).

Here is the main theorem we are aiming for.
 
\proclaim{Theorem 6.1}
Let $C$ be a $\Cal C^{k+\gamma}$ hyperbolic Cantor set, and let
$\Phi^y:I\to I$
be defined as in \S 5. We claim:
\item{(i)} if $k=1,2,\dots$, then $\Phi^y$ is $\Cal C^{k+\gamma}$
(with same H\"older exponent, but a different H\"older constant).
\item{(ii)} if $k=\infty $ or $\omega$, then
$\Phi^y$ is $\Cal C^{\infty}$ or $\Cal C^{\omega}$ respectively.
\endproclaim
 
First we need a few lemmas.
 
\proclaim{Lemma 6.2} For $A\subseteq \r$, if $f:A\to\r$ and $g:A\to\r$
are bounded, $\gamma$-H\"older continuous with H\"older constants $c,d$, then:
\item{(i)} $f+g$ is $\gamma$-H\"older with constant $c+d$, and
\item{(ii)} $f\cdot g$ is $\gamma$-H\"older
with constant $c\|g\|_\infty+
d\|f\|_\infty$.
\endproclaim
\demo{Proof}
(i) is immediate. The argument for (ii) comes by imitating the proof of
Leibnitz' rule in the Calculus:
$$
\align
|f(x)g(x)-f(y)g(y)| &=
|[f(x)-f(y)]g(x)+f(y)[g(x)-g(y)]| \\
&\leq c|x-y|^\gamma\|g\|_\infty
+d|x-y|^\gamma\|f\|_\infty.
\endalign
$$
\phantom{xxxxxxxxxxxxxxxxxxxxxxxxxxxxxxxxxxxxxxxxxxxxxxxxxxxxxmmmmmmmmmm} \ 
\qed
 
\enddemo
 
\proclaim{Lemma 6.3} Fix $k\geq 1$. For $A, B\subseteq \r$, let
$f:A\to B$ and $g:B\to\r$ be such that
$D^k(f) $, $D^k(g) $ are bounded and
$\gamma$-H\"older. Then $D^k(g\circ f) $ is $\gamma$-H\"older.
\endproclaim
\demo{Proof}
This now follows by induction, from the Chain Rule plus Lemma 6.2.
\ \ \qed \enddemo
 
The next lemma is basically the same as the ``bounded variation'' lemma, Lemma
1.15 from \cite{Bo 1}, except it is written in the reverse direction, for the
contractions $\varphi_i$ instead of the inverse map $S$. For the special case
$f_i=\varphi_{w_j}$
and $h_i=\log|D\varphi_{w_j}|$ where $i+j=n$,
one gets exactly the Bounded Distortion Property
(Theorem 2.1). The formulation given here is from the Appendix of \cite{Su 1};
the key idea for proving our Theorem 6.1, which is also in that Appendix 
(the sentence immediately preceding the Corollary there),
 will be how to
use this lemma to control higher order derivatives of the composition. One can 
summarize the idea as follows:
do not look at $\log|D^k\varphi_{x_i}|$, but instead at
$D^{k-1}\log|D\varphi_{x_i}|$. Then we are applying the linear operator $D^k$ 
to
a sum, which leads to the proof. 
 
As usual, for notational simplicity, we assume strict hyperbolicity.

\proclaim{Lemma 6.4} Consider a composition
of contractions $f_n\circ \dots\circ f_1$,
with $f_i :J_i\to J_{i+1}$
for intervals $J_i\subseteq \r$
, and such that $|Df_i|<\beta<1$.
For a point $x\in J$, write $x_1=x,\ x_{i+1}=f_i(x_i)$.
Let $h_i :J_i\to\r$
be H\"older continuous
with the same exponent $\gamma\in (0,1]$ and same constant $c$.
Then for $h(x)\equiv h_1(x_1)+\dots+h_n(x_n)$,
$h$ is also $\gamma$-H\"older continuous, with constant
$c_0=c\beta^\gamma/(1-\beta^\gamma)$ (independent of $n$).
\endproclaim
\demo{Proof}
Immediate from the geometric series, since for $x,y\in J_1$ we have
$|x_i-y_i|<\beta^i$.\ \ \qed \enddemo
 
The next little lemma is more subtle than one might at first think.
We wish to thank Z. Nitecki and M. Urbanski for discussions which
resulted in a first proof, and Y. Kifer for then finding the
much simpler argument given here.
 
\proclaim{Lemma 6.5}Let $f_n:I\to I$
be continuous functions with continuous $k^{\text{th}}$ derivative
and assume that there exist functions $f,g$ such that:
\item{(i)} $f_n\to f$ and
\item{(ii)} $D^k f_n\to g$, uniformly as $n\to \infty$.
Then $D^k f=g$.
\endproclaim
\demo{Proof}
We define, for each $0\leq j\leq k $, functions $g_j$,
and sequences of functions
$f_{n,j} $ and
$p_{n,j}$
by:
$$
\align
& g_k=g\text{ and } g_{j-1}(t)=\int_0^tg_j; \\
& f_{n,k}=D^kf_n
\text{ and } f_{n,j-1}(t)=\int_0^t f_{n,j}; \\
&  p_{n,j}=D^j f_n-f_{n,j}.
\endalign
$$
Thus
$$\align
& p_{n,k}\equiv 0, \\
& p_{n,k-1}(t)j=(D^{k-1} f_n)(t)-((D^{k-1} f_n)(t)-(D^{k-1} f_n)(0))\equiv
(D^{k-1} f_n)(0), \\
& p_{n,k-2}(t)=D^{k-2}f_n (0)+t(D^{k-1} f_n)(0)
\endalign
$$
and similarly,
(for each $n$) $p_{n,j}$
is for all $j$ a polynomial of degree $k-j-1$,
such that
$D p_{n,j}=p_{n,j+1}.$
Now for each $j$,
$\lim_{n\to \infty} f_{n,j}=g_j.$
In particular,
$$
g_0=\lim_{n\to \infty} f_{n,0}=\lim_{n\to \infty} (f_{n}-p_{n,0})
=f-\lim p_{n,0}.
$$
Hence $\lim_{n\to \infty} p_{n,0}$ converges (uniformly),
so to some polynomial $p_0$with degree at most $k-1$,
and we have
$g_0=f-p_0$.
Therefore,
$$
g=D^k g_0=D^k f-D^k p_0=D^k f
$$
as claimed.
\ \ \qed \enddemo
 
\demo{Proof of Theorem 6.1}
Since
$D^k \varphi_0$ and
$D^k \varphi_1$
are $\gamma$-H\"older,
by Lemma 6.3 so is
\newline
$D^{k-1}\log D \varphi_i$,
with some H\"older constant $c_0$.
Now we apply Lemma 6.4 to
$$
f_n\circ\dots\circ f_1=\varphi_{y_n}
\circ\dots\circ \varphi_{y_1}\equiv\varphi_{n}^y
$$ and
$$
h_j\equiv D^{(k-1)}\log D\varphi_{y_j}.
$$
For $x\in I$,
writing $x_1=x, \ x_2=h_1(x)$
etcetera as in Lemma 6.4, since
$$
D^{(k-1)}\log D\varphi_{n}^y=\sum_{j=1}^n h_j(x_j),
$$
we conclude that
$
D^{(k-1)}\log D\varphi_{n}^y$
is $\gamma$-H\"older,
with some different constant $c_1$ which is however independent of $n$.
 
Now to prove the Theorem, first  consider the case $k\geq 2$.
Here we have
$D^{(k-1)}\log D\Phi_{n}^y=D^{(k-1)}\log D\varphi_{n}^y$,
since the constant derivative of
$A^y_n$
disappears upon higher differentiation.
For $k=1$ these are not equal; they differ by the constant
$\log D A^y_n$
(which increases with $n$). However these cancel upon
subtraction, so in either case we have, for any $a,b\in I$,
$$
|\log D\Phi_{n}^y (a)-\log D\Phi_{n}^y (b)|=
|\log D\varphi_{n}^y (a)-\log D\varphi_{n}^y (b)|
.
$$
Therefore for all $k\geq 1$,
$D^{k-1}\log D\Phi_{n}^y$
is a sequence of bounded functions
which is $\gamma$-H\"older with the same constant, $c_1$.
Also, this sequence is uniformly bounded. For
$k=1$ this follows from bounded distortion, as in the proof of Theorem 5.4,
and in fact a bound is $e^{K\beta^\gamma}$.
For $k>1$,
we argue as follows:
if it were unbounded, then by H\"older continuity with the same constant,
some subsequence goes uniformly to either $+\infty $ or
$-\infty$. By integration $(k-1)$ times, by induction
 this contradicts the boundedness for
$k=1$.
This implies equicontinuity. Now by boundedness and equicontinuity,
 there  is some convergent
subsequence, using the standard diagonalization argument as in the proof of the
Arzela-Ascoli Theorem.
At the same time, from \S 5 we know that
$\lim  D\Phi_{n}^y = D\Phi^y $
exists which implies
$\log D\Phi_{n}^y $
converges to
$\log D\Phi^y.$
Calling the subsequence
$$
\log D\Phi_{n_j}^y=f_j,
$$
we are in the situation of Lemma 6.5:
$f_j\to f$,
$D^{k-1}f_j\to g$
hence
$D^{k-1}f= g$.
Thus
$D^{k-1}\log D\Phi^y$
is a uniform limit of $\gamma$-H\"older functions with the
same constant $c_o$,
hence the limit is $\gamma$-H\"older with
constant $c_o$.
From Lemma 6.3,
$
D^{k}\Phi^y$
is also
$\gamma$-H\"older and we are done for $k= 1,2,\dots$.
 
Finally note that for $k=\infty$
we are done by part (i),
and for $k=\omega$
we can apply Arzela-Ascoli to see that
$\{ \Phi_n^y \}$
is a normal family,
hence the limit
$ \Phi^y $
is also analytic.
\ \ \qed \enddemo

\subhead{Remark:}\endsubhead  
We emphasize again the subtle point in the logic of this argument:
$\Cal C^{1+\gamma}$ 
convergence of $ \Phi_n^y $ to $ \Phi^y $
is \it not \rm known. What we \it do \rm know is convergence in the 
$\Cal C^1$ norm (from Theorem 5.9) and convergence along a \it subsequence \rm
in the $\Cal C^{1+\gamma}$ norm, as just shown. This is enough to prove 
the claim of the Theorem.

\subheading{\S 7. Smoothness of limit sets, and rigidity}
 
Given two $\Cal C^{1+\gamma}$ 
hyperbolic Cantor sets $(C,S)$
and $(\widehat C,\widehat S)$,
recall that
 the (full) conjugacy $\Phi$ is an order-preserving map
defined on all of $I$.
This map is uniquely determined on $C$
by the conjugacy equation,
since, as one sees, the symbolic dynamics is preserved.
Note that for any two topological Cantor sets,
once they have been coded by the two-shift $\Sigma^+$
in an order-preserving way, this conjugacy on the Cantor sets extends
to a homeomorphism on $I$.
The issue therefore is what types of conjugacies preserve what type of 
structure.
As is well known
and not hard to show, for instance, a biLipschitz
$\Phi$ will preserve the Hausdorff
dimension. We noted in \cite{BF 1} that $\Cal C^1$ maps preserve the order-two 
density.
Furthermore for $\Cal C^1$ conjugacy from \cite{Su 1} one has:

\proclaim{Lemma 7.1}
   If two $\Cal C^{1+\gamma}$ hyperbolic Cantor sets are
$\Cal C^1$ conjugate, then they have the same scaling function.
\endproclaim
\demo{Proof}
By uniform continuity of the derivatives, since we already know 
the scaling functions exist from Theorem 3.1,
this is immediate.
\ \ \qed
\enddemo 
 
Hence under the same assumption, by Corollary 5.3 we have:
\proclaim{Corollary 7.2} They have the same collection of limit sets.
\ \ \qed\endproclaim

To prove our rigidity theorem, we will need the following.
 
\proclaim{Lemma 7.3} Let $(C,S)$
be a  hyperbolic $\Cal C^{k+\gamma}$
Cantor set.
Let $\widehat S:I_0\cup I_1 \to I$
be a $\Cal C^{k+\gamma}$
map such that
$S=\widehat S$ on $C$. Then
$S$ and $\widehat S$ are
$\Cal C^{k+\gamma}$ conjugate.
\endproclaim
\demo{Proof}
The conjugacy is the identity
map on $C$; what we want to do is define it on the gaps.
We begin by defining
$\Phi$ to be the identity also on the gap $G$ between
$I_0 $ and $I_1$.
The conjugacy is then uniquely defined from the conjugacy equation,
by the dynamics.
That is, writing
$G_{x_0\dots x_n}=\varphi_{x_0\dots x_n}(G)$,
we have for
$a\in G_{x_0\dots x_n}$,
$$
\Phi(a)=\widehat\varphi_{x_0\dots x_n}
(\varphi^{-1}_{x_0\dots x_n}(a))=
\widehat\varphi_{x_0\dots x_n}(S^n(a)).
$$
One immediately checks that with this definition,
$\Phi$ is a conjugacy.
 
This map is $\Cal C^{k+\gamma}$ on the interiors of all the gaps.
At points in $C$, to check $\Cal C^{k+\gamma}$ one must be careful because 
these
points  
are also limits of interior points in the gaps.

Here is one way of proving $\Phi$
is everywhere $\Cal C^{k+\gamma}$.
 
Define a sequence of maps
$\Phi_n:I\to I$
by:
$\Phi_0=$ the identity,
$\Phi_1=\Phi_0$ on $G$ and $\widehat\varphi_{x_0}\circ S$
everywhere else (i.e\. on $I_0\cup I_1$),
and inductively, set
$\Phi_n$ to be equal to $\Phi_{n-1}$
everywhere except on
$\bigcup I_{x_0\dots x_n}$,
where it is defined to be
$\widehat\varphi_{x_0\dots x_n}\circ S^n$.

These maps converge uniformly to $\Phi$.
So if we can show
that for each $n$, $D^k \Phi_n$
is $\gamma$-H\"older with a constant independent of n, this will carry over to 
the limit and we 
will be done.
(Here we will use the fact that the maps
$\widehat\varphi_{x_0\dots x_n}\circ S^n$
are
$\gamma$-H\"older with a fixed constant).

The advantage of this method is that we must only check smoothness at each 
stage,
and so each time at only finitely many points.
 
Now consider the map $f=\widehat\varphi_{x_0\dots x_n}\circ S^n$
on $C\cap I_{x_0\dots x_n}$.
It is the identity there, and since $C$ is dense in itself,
$Df=1$ on that set.
Since it is twice differentiable,
$D^2f=0$ there and similarly for $D^kf$.
Therefore when we define $\Phi_n$
by gluing together
$\widehat\varphi_{x_0,\dots,x_{n-1}}\circ S^{n-1}$ and
$
\widehat
\varphi_{x_0\dots x_n}\circ S^n$
at an endpoint $p$,
the two functions agree at $p$ for all derivatives $\leq k$.
Also $D^k$ is $\gamma$-H\"older, for each piece.
Hence for all $n$,
$\Phi_n$ is $\Cal C^{k+\gamma}$ with a fixed H\"older constant, as we 
wanted to show; so we are done.\ \ \qed\enddemo

This produces one conjugacy. 
In \S 8 we  will return to this proof, in order to study 
\it how many \rm such maps $\Phi$ there are.

We are now ready to prove:
 
\proclaim{Theorem 7.4 (highest smoothness)}
Given a $\Cal C^{1+\gamma}$ hyperbolic Cantor set
$C$,
its limit sets have the highest degree of smoothness
of any hyperbolic
$\Cal C^{1+\gamma}$ Cantor set in the $\Cal C^{1+\gamma}$- conjugacy class of $C$.
\endproclaim
\demo{Proof}
Let $(\widehat C,\widehat S)$
be a
 $\Cal C^{k+\gamma}$,
$\Cal C^\infty$ or $\Cal C^\omega$
hyperbolic Cantor set
which is $\Cal C^1$ conjugate to $C$.
By Corollary 7.2, $C$ and $\widehat C$
have the same limit sets.
And by Theorem 6.1, the map
$\widehat \Phi^y:I\to I$ defines a dynamics
$S^y:I^y_0\cup I^y_1\to I$
by conjugation with the map $\widehat S$,
which has that same degree of smoothness.
\ \ \qed\enddemo

\proclaim{Theorem 7.5 (rigidity)}
If $(C,S)$
and $(\widehat C,\widehat S)$ are two
$\Cal C^{k+\gamma}$,
$\Cal C^\infty$ or $\Cal C^\omega$
hyperbolic Cantor sets which
either (a) are $\Cal C^1$
conjugate by a map $\Phi$,  or
(b) have the same scaling function $R$,
then they are in fact conjugate by a map
$\widetilde \Phi:I\to I$
which for (a) agrees with $\Phi$ on $C$ or 
for (b) agrees with the coding;
this map is
$\Cal C^{k+\gamma}$,
$\Cal C^\infty$ or $\Cal C^\omega$ respectively.
\endproclaim
\demo{Proof}
By either hypothesis they have the same limit sets. Choose one, $C^y$.
Again by Theorem 6.1,
the maps
$\Phi^y$, $\widehat\Phi^y$
have the same smoothness as $S$, $\widehat S$.
Now let
$S^y$, $\widehat S^y$
denote the maps defined on
$I^y_0\cup I^y_1$ by these conjugacies.
We are exactly in the situation of Lemma 7.3, and
have a conjugacy $\Phi$
of
$S^y$ and $\widehat S^y$.
Composing the three maps
$$
(\widehat \Phi^y)^{-1}\circ\Phi\circ \Phi^y
$$
finishes the proof.
\ \ \qed\enddemo

\subheading{\S 8 Banach space structure}

Fix a hyperbolic $\Cal C^{1+\gamma}$
Cantor set $(C,S)$. For 
$r=k+\gamma $ where $\gamma\in (0,1],\ k\geq 1$ or for $r=\infty,\ \omega$
we write $\Cal E^r\equiv \Cal E^r(C)$ for the collection of Cantor sets (with
maps)
which are $\Cal C^r$- conjugate to $(C,S)$. (From Lemma 1.1.2, these are also 
hyperbolic $\Cal C^{1+\gamma}$ Cantor sets). 
The spaces
$\text{Diff}^{r}$,
$\wh{\Cal E}^{r},$
$\Cal E^{r}_*$ and $\Cal E^{r}_{**}$  are defined as they were in
\S 4 for the case $r= 1+\gamma$. 

In this section we will see how $\wh{\Cal E}^{r}$
 can be viewed as a Banach manifold, 
in fact a Banach Lie group. We will also define a natural
topology on ${\Cal E}^{r}$, and show that $\wh{\Cal E}^{r}$
factors nicely over ${\Cal E}^{r}$ as a topological space.

We define first the $\Cal C^\gamma$ norm on the H\"older functions
$\Cal C^\gamma(I,\r)$ to be
$$
\|f\|
_{\Cal C\gamma}
=\|f\|_\infty
+\sup_{x,y\in I}\frac{|f(x)-f(y)|}{|x-y|^\gamma}.
$$
For  $r=k+\gamma $ where $\gamma\in (0,1],\ k\geq 1$, the  $\Cal C^r$
norm will be 
$$
\Sigma_{l=0}^{k-1} \|D^l f\|_\infty 
+ \|D^k f\|
_{\Cal C\gamma}.
$$
For $\Cal C^\infty$ we define 
$$\|f\|
_{\Cal C\infty}=\sup_l\{\|D^l f\|_\infty \},$$
and for $\Cal C^\omega$
we will use the sup norm (since it is  equivalent to all the other 
$\Cal C^r$ norms there).

By definition a Banach manifold is a manifold which is locally 
modelled on a Banach space, and a Lie group is a group, which is also a 
$C^\infty$ manifold modelled on a complete, locally convex vector space
(see e.g. \cite{Mi}).
Recall that $\text{Diff}^r$
denotes the $\Cal C^r$ order-preserving diffeomorphisms of $I$.
Now choice of a set in $\Cal E^r$ 
identifies the collection $\wh{\Cal E}^r$ 
with
$\text{Diff}^r$,
as we have seen in \S 4. $\text{Diff}^r$ 
is an open subset of 
$\Cal C^r_{0,1}(I,\r)$, 
which is how we will write 
the set of all $\Cal C^r$
functions from $I$ to $\r$ such that
$f(0)=0$ and $f(1)=1$.
This in turn is a closed affine subspace of 
$\Cal C^r(I,\r)$.
To see this note that, defining 
$$
\Cal B^r_{0,1}(I,\r)
=\{ f\in \Cal C^r(I,\r)
:\ f(0)=0=f(1)\},$$
two functions in 
$\Cal C^r_{0,1}$
differ exactly by an element of 
$\Cal B^r_{0,1}$.
Now  
$\Cal B^r_{0,1}$ is a Banach space,  with the $\Cal C^r$  norm.
Hence $\wh{\Cal E}^r$  is a Banach manifold:
it is identified with $\text{Diff}^r$,
which in turn corresponds to an open subset of 
$\Cal B^r_{0,1}$.
Now $\text{Diff}^r$ is a group,
hence it (and therefore $\wh{\Cal E}^r$) 
is a Banach Lie group.
Two choices have been made:
the choice of a Cantor set in $\Cal E^r$,
and of a special point (the identity) in  $\text{Diff}^r$.
These choices determined the maps to $\Cal B^r_{0,1}$
and hence the  metric (inherited from the $\Cal C^r$ norm).
Both choices moreover amount to the same thing:
changing $C$ to $D$ in $\Cal E^r$ 
(as  in Proposition 5.1, for $r=1$)
corresponds to
a right translation in the group $\text{Diff}^r$.

Now in a Lie group
one ideally would like to work with a (left- or right-) invariant metric.
If the group is compact (or, more generally, amenable),
one can make a given metric invariant (while keeping equivalence)
by averaging over translations.
In our case, however, one cannot get an equivalent 
invariant metric- $\text{Diff}^r$ 
is not only non-compact but non-amenable! The (non-uniform) bounded equivalence
proved in  Proposition 5.1 is nevertheless enough for what we needed,
for the proof of Theorem 5.3.

In summary we have:
\proclaim{Proposition 8.1} 
$\wh{\Cal E}^r$   
is a Banach manifold. It is 
naturally identified up to right composition with the Banach Lie group
$\text{Diff}^{\, r}$, 
and with an open subset of a closed affine subspace of 
$\Cal C^r(I,\r)$.
\ \ \qed
\endproclaim

A similar estimate to that shown in Proposition 5.1 for $r=1$ 
holds for $r>1$. Therefore one has, for $r=k+\gamma, \infty, \omega:$

\proclaim{Proposition 8.2} 
\item{(a)}
The $\Cal C^r$ metric 
on $\text{Diff}^{\ r}$ is right-invariant up to
(non-uniform) bounded equivalence.
\item{(b)}
The  $\Cal C^r$ metric on $\wh{\Cal E}^r$   is base-point independent 
up to (non-uniform) bounded equivalence.
\phantom{The  $\Cal C^r$ metric on $\Cal E^r$  is base-point independent 
up to (non-uniform) bounded equival}
\ \ \qed\endproclaim

We note that $\Cal E^r_*$, the space of marked Cantor sets, is also
a Banach 
manifold, 
by the same reasoning as for $\text{Diff}^{ r}$: 
it is an open  subset of a  closed affine subspace of
$\Cal C^r(C,\r)$.

Next we will describe more fully the relationship between the spaces 
$\wh{\Cal E}^r$  and ${\Cal E}^r$. 
For $r=k+\gamma$, 
we write  
$\text{Diff}^r_0 (I)$ 
for the  collection  of $\Cal C^r$ diffeomorphisms of the 
unit interval $I$
whose first $k$ derivatives are $1,0,\dots,0$  at the endpoints. 
This is also a Banach manifold.

\proclaim{Proposition 8.3} Given the choice of a Cantor set $C$,
$\wh{\Cal E}^r$ 
factors naturally, set theoretically and topologically, as
$$
\wh{\Cal E}^r=\Cal E^r\times \text{Diff}^r_0  (I),
$$ 
with the  topology on $\Cal E^r$  defined below.

\endproclaim
\demo{Proof}
First,
 let us consider 
how many $\Cal C^{r}$ 
maps from
$I$ to $I$ there are, conjugating   $(C,S)$
with $(C,\widehat S)$ (the Cantor sets are the same, but the maps
may be different off of the Cantor sets).
 In the proof of Lemma 7.3, note that instead of starting with 
$\Phi$ equal to the identity on the gap $G$ we could have taken
any $\Cal C^{k+\gamma}$ 
diffeomorphism from $G$ to itself, whose derivatives
agree
with the identity at the endpoints, up to order $k$.
Conversely any conjugacy is specified by its values on
$G$, since elsewhere it is then determined by the dynamics. 
Therefore we see that  the set of $\Cal C^{r}$ 
conjugacies
from   $(C,S)$
to $(C,\widehat S)$ correspond naturally to $\text{Diff}^r_0 (I)$.
(This is also true when $S=\wh S$!)

Next we
 consider  how many 
different $\Cal C^{r}$ conjugacies are possible 
in the rigidity theorem, from $(C,S)$ to $(\widehat C,\wh S)$.
 Given the existence of  one such 
map and hence a restricted 
conjugacy, we can define (all the) 
other extensions by a method like that used in 
the proof just given. That is, we first define the 
conjugacy arbitrarily on the first-level gaps (but with the correct 
derivatives of order $\leq k$ at the endpoints). Then we extend 
by the dynamics.
Or, we can quote that statement directly,
making use of use a ratio Cantor set as intermediary as in the proof of 
Theorem
7.5, and now replacing $\Phi$ by one of the more general maps
 described above.

This shows we have a product of sets. 
$\Cal E^r$ 
has not yet been given a topology.
But from the product decomposition, we can define a 
family of metrics as follows.
Choosing one element 
of 
$\text{Diff}^r_0 $ 
defines an embedding into 
$\wh{\Cal E}^r$, and we just use the $\Cal C^r$ metric there.
(One would like to get a more natural 
definition by  taking the infimum  over all such choices; however, 
 it is then not clear that the triangle inequality will  hold).
At any rate the metrics are equivalent, so this defines a 
natural topology on $\Cal E^r$.

We will show that the $\Cal C^r$ metric on $\wh{\Cal E}^r$
is equivalent to the product of the metrics on 
$\Cal E^r$ and 
$\text{Diff}^r_0 $.

It is easy to see that the map from  $\Cal E^r$
to each factor is continuous. (To $\text{Diff}^r_0 $ it is also affine).
For the converse, given the base point $(C,S)$,
let first $f,g\in \text{Diff}^r$ 
be such that 
$(C_f,S_f)=(C_g,S_g)$.
Write $f_0,\ g_0$ 
for the corresponding elements  of $\text{Diff}^r_0 $, i.e\.
the restrictions of $f$ and $g$ to the middle gap of $C$ (rescaled in the range).
We claim that if $f_0$ and $g_0$ 
are close in $\text{Diff}^r_0 $
then $f$ and $g$ 
are close in $\text{Diff}^r$.
The formula for $f$ on an $n^{\text{th}}$ level gap
of $C$ is
$$
f(a)= 
\varphi_{x_0\dots x_n}^f\circ f_0\circ
\varphi^{-1}_{x_0\dots x_n}(a).
$$
Here $\varphi,\ \varphi^f$
denote inverse branches for
$S$ and $S_f$ respectively.
By assumption $\varphi^f=\varphi^g$.
Now by bounded distortion (Lemma 6.4) for $k=1$, and for general $k$ by 
the proof of Theorem 6.1,
$D^k \varphi_{x_0\dots x_n}^f$ 
is uniformly $\gamma-$H\"older
with constant independent of $n$.
This proves $\|f-g\|_{\Cal C 1}$
is small, which is what we wanted to show.

Next we drop the assumption that $f$ and $g$ 
give the same maps.
We have chosen an element of $\text{Diff}^r_0 $
to define the metric on $\Cal E^r$.
 Let $\tilde f,\ \tilde g$ 
denote the maps in $\text{Diff}^r$
such that $\tilde f_0=\tilde g_0$ 
is that element, with $(C_f,S_f)=(C_{\tilde f},S_{\tilde f})$
and similarly for $g$.
Now by definition, the distance between the pairs in
$\Cal E^r$ is $\|\tilde f-\tilde g\|_{\Cal C 1}$.
So we just apply the triangle inequality using the previous case,
to conclude that
$d_C(C_f,C_g)\equiv \|f-g\|_{\Cal C 1}$ is also small.
\ \ \qed
\enddemo

\subheading{Remark} 
In a conversation about the proof of Lemma 7.3,
 Yair Minsky  pointed out to us an interesting parallel 
between 
that argument and 
Sullivan's ``flexibility and rigidity'' theorem for 
Kleinian groups. 
Sullivan  showed, for a finitely generated Kleinian 
group $\Gamma$,
that the limit set $\Lambda$ of the group 
itself is ``rigid'',  i.e\. 
a 
quasiconformal conjugacy (to another Kleinian group) which lives
(Lebesgue almost-surely) on $\Lambda$ must be M\"obius. This is 
a consequence of Sullivan's lemma that $\Lambda$   carries  no measurable
$\Gamma$-
invariant line fields. 
(Note that by contrast, for hyperbolic
Cantor sets, quasisymmetric conjugacy does not imply
smooth conjugacy; as we have seen, 
one also needs  to know the scaling function). 

 Sullivan used this to show that a quasiconformal conjugacy
is determined by a Beltrami differential on $\Omega / \Gamma$ 
where 
$\Omega$ is the domain of discontinuity. Thus, in a sense, one  
has rigidity on the limit set and flexibility off of it.
The group actions correspond to the
two
(\it not \rm restricted) expanding maps, and the surface
$\Omega / \Gamma$, or equivalently a fundamental domain for the 
action of $\Gamma$ on $\Omega$,  is analogous
to the gap $G$ of the Cantor set. As is the case there,
 the conjugacy is then 
specified elsewhere 
by the dynamics. Sullivan's theorem can then be stated as follows: 
Teich($\Gamma$)$=$ Teich($\Omega / \Gamma$), where this refers to the
Teichm\"uller space of a group and of a surface respectively; 
 this
formulation led to the statement in the Proposition above.

\subhead
{Concluding remarks: limit sets as the attractor of a semigroup action}
\endsubhead 
$\Cal E^{1+\gamma}$ denotes  the $\Cal C^{1+\gamma}$ equivalence 
class of a given 
$\Cal C^{1+\gamma}$ hyperbolic
Cantor set. 
 The nested subclasses $\Cal E^r$,
for maps of smoothness $r=k+\gamma,\  \infty,\  \omega$, 
also are conjugate 
with that higher degree of smoothness. 
Thus smoothness classes are also conjugacy classes. Choosing one set in 
$\Cal E^r$ as a base point,  $\Cal E^r$ is naturally identified with 
a topological factor of 
the $\Cal C^r$ orientation-preserving
diffeomorphisms of the interval,
$\text{Diff}^r$
which is a Banach 
manifold. Moreover we can choose one Cantor set as a common base point
for all the $\Cal E^r$, since by Theorem 7.4 smoothest Cantor sets exist.
Then the nested collections $\Cal E^{1+\gamma}
\supseteq \dots \supseteq \Cal E^r\dots$
are naturally identified with factors of
$\text{Diff}^{1+\gamma}
\supseteq \dots \supseteq 
\text{Diff}^r\dots$. 
(Each is a Banach manifold with its own topology,
and is a dense subset of the larger collections, with respect to their
topologies).
The  spaces of marked Cantor sets $\Cal E^r_*$ are also Banach manifolds.
The free semigroup
on two generators  $FS_2$ acts on each submanifold $\Cal E^r_*$
by replacing  it with its left or right Cantor subset. 
From Theorem 7.4, the limit sets are in the intersection of the 
$\Cal E^r_*$.
  From Theorem 5.3, because the bounds are uniform over all 
Cantor subsets of level $n$,
 the collection of limit sets is 
an attractor for this action. This convergence is
exponentially fast in the $\Cal C^1$ norm. (Warning: we  have only shown
convergence in \it this \rm norm; see the Remark at the end of \S6).
The semigroup action on the attractor itself
can be described symbolically very simply as follows. Recall
 the map $y\mapsto C^y$ for $y$ in
 the dual Cantor set 
$\Sigma^-$ and $C^y$ the corresponding ratio Cantor set. Now just
concatenate  $y$ on the right
with a finite string of symbols.

\enddocument